\documentclass[11pt]{article}
\usepackage[dvips]{color}
\usepackage{hyperref}
\usepackage{mathrsfs}
\usepackage{graphicx}
\usepackage{amsfonts}
\usepackage{amssymb}
 \usepackage{epstopdf}
 \usepackage{hyperref}
\usepackage{multirow}
\usepackage[pagewise]{lineno}
\usepackage{tikz}
\usetikzlibrary{arrows,shapes,chains}

\definecolor{dgreen}{rgb}{0,.8,.3}
\definecolor{blue}{rgb}{.2,.3,.7}
\definecolor{red}{rgb}{1.0,0.2,0.2}
\usepackage{color}
\input amssym.def

\input epsf.tex

\begin{document}

\renewcommand{\Box}{\rule{2.2mm}{2.2mm}}
\newcommand{\BOX}{\hfill \Box}

\newtheorem{eg}{Example}[section]
\newtheorem{thm}{Theorem}[section]
\newtheorem{lemma}{Lemma}[section]
\newtheorem{example}{Example}[section]
\newtheorem{remark}{Remark}[section]
\newtheorem{proposition}{Proposition}[section]
\newtheorem{corollary}{Corollary}[section]
\newtheorem{defn}{Definition}[section]
\newtheorem{alg}{Algorithm}[section]
\newtheorem{ass}{Assumption}[section]
\newenvironment{case}
    {\left\{\def\arraystretch{1.2}\hskip-\arraycolsep \array{l@{\quad}l}}
    {\endarray\hskip-\arraycolsep\right.}

\def\argmin{\mathop{\rm argmin}}

\makeatletter
\renewcommand{\theequation}{\thesection.\arabic{equation}}
\@addtoreset{equation}{section} \makeatother

\title{Relaxed constant positive linear dependence constraint qualification for disjunctive systems}
\author{ Mengwei Xu\thanks{\baselineskip 9pt Institute of Mathematics, Hebei University of Technology, Tianjin 300401, China. 
E-mail: xumengw@hotmail.com. This author is supported by the NSFC  No. 11901556, No. 12071342 and the HNSF No. A2020202030} \ and \   Jane J. Ye\thanks{\baselineskip 9pt Corresponding author, Department of Mathematics
and Statistics, University of Victoria, Victoria, B.C., Canada V8W 2Y2. E-mail: janeye@uvic.ca.
This author is partially supported by NSERC.}
}
\date{}
\maketitle
{\bf Abstract.} The disjunctive system is a system  involving a disjunctive set which is  the union of finitely many polyhedral convex sets. In this paper, we introduce a notion of  the relaxed constant positive linear dependence constraint qualification (RCPLD) for  the disjunctive system.  {For a disjunctive system,  our notion  is weaker than the one we introduced for a more general system}
recently  (J. Glob. Optim. 2020) and is still  a constraint qualification.
To obtain the local error bound for the disjunctive system, we introduce the  piecewise RCPLD under which the error bound property holds if all inequality constraint functions are subdifferentially regular and the rest of the constraint functions are smooth. We then specialize our results to the ortho-disjunctive program, which includes the  mathematical program with equilibrium constraints (MPEC), the mathematical program with vanishing constraints (MPVC) and the mathematical program with switching constraints (MPSC) as special cases. For MPEC,  we recover MPEC-RCPLD, an MPEC variant of RCPLD and propose the MPEC piecewise RCPLD to obtain the {error bound property}. For MPVC,  we introduce new constraint qualifications MPVC-RCPLD and the piecewise RCPLD, which also implies  the local error bound. For MPSC, we show that both RCPLD and the piecewise RCPLD coincide and hence it leads to
 the local error bound.



\newpage

\baselineskip 18pt
\parskip 2pt

\section{Introduction}
This paper considers the  disjunctive system of the form
\begin{eqnarray}\label{feasibility}
 g(x)\leq 0,
 h(x)=0,
\Phi_i(x)\in {\Gamma_i}, & i=1,\cdots,l,
\end{eqnarray}
where  $g:\mathbb{R}^d\to \mathbb{R}^n$ is locally Lipschitz continuous,   $h:\mathbb{R}^d\to \mathbb{R}^m$, $\Phi_i(x): \mathbb{R}^d\to \mathbb{R}^{p_i}$ are  smooth  and each $\Gamma_i \subseteq \mathbb{R}^{p_i}$
 is a disjunctive set, i.e., 
the union of finitely many polyhedral convex sets, $i=1,\cdots, l$. 
Let ${\cal F}$ be  the solution set of the disjunctive system $(\ref{feasibility})$. 

The mathematical program with disjunctive constraints  (MPDC) is an optimization problem with the constraint region defined by  the disjunctive system (\ref{feasibility}). 
Several classes of mathematical programs can be reformulated as MPDC, including some mixed-integer programs (\cite{balas}) and the so-called ortho-disjunctive programs recently introduced in Benko et al$.$ \cite{bch}, where the disjunctive set is the union of finitely many boxes. Several classes of important optimization problems can be reformulated as  ortho-disjunctive programs.
 Indeed, when 
$$\Gamma_i:=\Omega_{E}:=\left \{(y,z)\in \mathbb{R}^2: yz=0, y\geq 0,z\geq 0\right \}=\mathbb{R}_+\times \{0\} \cup \{0\} \times\mathbb{R}_+,$$ the disjunctive system (\ref{feasibility}) becomes the equilibrium/complementarity system and MPDC becomes the mathematical program with equilibrium constraints (MPEC)  (\cite{lpr,Outra}); when 
$$\Gamma_i:=\Omega_{V}:=\left \{(y,z)\in \mathbb{R}^2: z\geq 0, yz\leq 0 \right \}=\mathbb{R}_-\times \mathbb{R}_+\cup {\mathbb{R}} \times  \{0\},$$ the disjunctive system (\ref{feasibility}) becomes the vanishing system and MPDC becomes 
 the mathematical programs with vanishing constraints (MPVC) (\cite{ak,hk08,hk09,hks12,kn}) and when 
$$\Gamma_i:=\Omega_{S}:=\left \{(y,z)\in \mathbb{R}^2: yz= 0 \right \}=\mathbb{R}\times \{0\} \cup \{0\}\times\mathbb{R},$$ the disjunctive system (\ref{feasibility}) becomes the switching system and MPDC is
 the mathematical program with switching constraints (MPSC) (\cite{kms,lg,ly,mmpsc}).
 

Necessary optimality conditions for  MPEC, MPVC and MPSC are notoriously difficult to obtain 
when they are treated as nonlinear programming problems (see e.g.  \cite{ak,mmpsc,yzz}). On the other hand, tailored constraint qualifications may hold when they are treated as MPDC.
Hence one way to copy with this difficulty is to rewrite these problems as MPDCs, apply constraint qualifications for the corresponding MPDCs and hence obtain tailored constraint qualifications
 for these problems. So far there are many  MPDC variants of  constraint qualifications and optimality conditions have been introduced and studied  (see e.g. \cite{bch,bg, fko, g14, mam, mmpscnew, tca}); see \cite[Section 3]{ly} for a survey on most of these conditions.

The relaxed constant positive linear dependence (RCPLD) condition was  first proposed by Andreani et al.  \cite[Definition 4]{RCPLD} as a constraint qualification for classical nonlinear programs.  RCPLD is  very useful for smooth nonlinear programs since it is weaker than the Mangasarian-Fromovitz constraint qualification (MFCQ) and leads to the existence of a  local error bound (\cite[Theorem 4.2]{gzl}).
 RCPLD was extended to a smooth system involving an abstract set by Guo and Ye \cite[Corollary 3]{GuoYe} and then weakened by
Xu and Ye  \cite[Definition 1.1]{xy2} to a nonsmooth system  with an arbitrary abstract set: $g(x)\leq 0, h(x)=0, x\in \Gamma$,  $g$ is locally Lipschitz continuous, $h$ is smooth and $\Gamma$ is closed. In \cite[Theorem 3.1]{xy2},   RCPLD for such a general nonsmooth system was proved to be a constraint qualification for the Mordukhovich (M-) stationarity condition.
%
Moreover, when  $g$ is subdifferentially regular/Clarke regular (which holds  e.g.  if either $g$ is  convex or smooth) and $\Gamma$ is subdifferentially regular/Clarke regular (which holds e.g.  if $\Gamma$ is  convex),
 it was shown in \cite[Theorem 3.2]{xy2} that  RCPLD for the nonsmooth system guarantees the existence of a local error bound.

A motivation for studying RCPLD for  disjunctive system (\ref{feasibility}) with nonsmooth problem data is to provide a constraint qualification for certain MPDCs which may not satisfy the usual constraint qualifications such as the no nonzero abnormal multiplier constraint qualification (NNAMCQ) (see e.g.,  \cite{byz}), which is also referred to the Mordukhovich criterion \cite{m1} and is equivalent to MFCQ in the case of nonlinear programming. 
One of such applications is  the bilevel programming problem for which it is well-known that NNAMCQ/MFCQ never hold and hence it is important to find other constraint qualifications for it. Recently, a few point-based constraint qualifications, which are weaker than NNAMCQ  have been introduced;  see e.g., the first  order sufficient condition for metric subregularity (FOSCMS) (see e.g., \cite[Theorem 4.3]{g13}),  quasinormality and the directional quasinormality (\cite{byz,bch}). However it was shown that even FOSCMS fails to hold  \cite[Proposition 5.1]{by} while it is possible for RCPLD to hold \cite{xy2} for the bilevel programming problems.
 Another motivation to study RCPLD for the disjunctive program is to derive some RCPLD-type conditions for local error bound. Although RCPLD leads to the existence of a local  error bound  for a smooth nonlinear program, there is no proof that it is still a sufficient condition when a non-subdifferentially regular abstract set is included in the program. In particular, since a disjunctive set  is  in general not Clarke regular, RCPLD for the disjunctive program may not imply the existence of a local  error bound. Although we do not have a proof for this conjecture, we notice that for MPEC, all RCPLD-type {constraint qualifications} in the literature require extra conditions on top of  RCPLD; see e.g., \cite[Theorem 5.1]{gzl} and  \cite[Theorem 3.2]{cl14}.  

In this paper we will  take the special structure of the disjunctive set into
account to provide a weaker and more precise notion of RCPLD than the one in Xu and Ye  \cite{xy2} designed for a nonsmooth system with an arbitrary set. To motivate our approach, consider a simpler system
$$\Phi(x)\in C:=\{y: \langle c_j, y \rangle \leq \alpha_j, j\in {\cal I}\}.$$ Obviously this system can be described by a system of inequalities, i.e.,
$$\Phi(x)\in C \Longleftrightarrow \langle c_j, \Phi(x) \rangle \leq \alpha_j,\quad  j\in {\cal I}.$$ Let the active set at $y\in C$ be ${\cal I}(y):= \{j\in {\cal I}:\langle c_j, y\rangle = \alpha_j\}$. By RCPLD for the system of inequalities,  a feasible solution $\bar x$ satisfies RCPLD if 
for any $I\subseteq {\cal I}(\Phi(\bar x))$ and
$\lambda_j \geq 0\ (j\in I)$ not all zero such that
\begin{equation}\label{eqn1.2}
0=\sum_{j\in I} \lambda_j \nabla \Phi(\bar x)^T c_j,\end{equation}
the set of vectors $\{\nabla \Phi( x)^T c_j\}_{j\in I}$ is linearly dependent for every $x\in U(\bar x)$, a neighborhood of $\bar x$. 
Since $C$ is a convex polyhedral set,  by \cite[Theorem 6.46]{var}, the normal cone to $C$ at $y$ is equal to
$${\cal N}_C (y)=\left \{ \sum_{j\in {\cal I}(y)} \lambda_j c_j: \lambda_j \geq 0, j\in {\cal I}(y)\right \}=cone(A_C(y)),$$ 
where 
$A_C(y):= \{c_j\}_{j\in \widetilde{\cal I}(y)}$ with $
\widetilde{\cal I}(y)\subseteq {\cal I}(y).$ 
Hence (\ref{eqn1.2}) can be rewritten as
\begin{equation}\label{eqn1.3}
0=\nabla \Phi(\bar x)^T \bar \eta\end{equation}with $\bar\eta=\sum_{j\in I} \lambda_j c_j \in {\cal N}_C (\Phi(\bar x))$ and  RCPLD at $\bar x$ is equivalent to saying that if there is a nonzero vector $0\not =\bar \eta \in {\cal N}_C (\Phi(\bar x))$ such that (\ref{eqn1.3}) holds, then the set 
$$\{\nabla \Phi( x)^T c_j\}_{j\in I}=\left \{\nabla \Phi( x)^T \beta \right \}_{\beta\in A} ,$$
where $A=\{ c_j: j\in I\}\subseteq A_C(\Phi( x))$ satisfying $\bar \eta \in {cone(A)}\subseteq {\cal N}_C(\Phi(\bar x))$,  must be linearly dependent. This is exactly the concept of RCPLD we introduced in this paper in Definition \ref{nrcpld} for this simplier system. Although  $\Gamma_i$ in the general system (\ref{feasibility}) is not convex polyhedral in general, the limiting normal cone is still generated in certain sense by finitely many vectors. This is the generator approach we will take in this paper.

Moreover in order to study error bounds for the disjunctive system, we will adopt 
the following local decomposition approach. For easy explanation, assume that 
  $\bar x$ is a solution to the system (\ref{feasibility}) with $l=1$ and $\Phi_1(\bar x) \in \Gamma_1=\bigcup_{r=1}^R C_r,$ where each $C_r$ is  convex polyhedral. Then  $\bar x$ is also a solution to each subsystem 
  $$g(x)\leq 0, h(x)=0, \Phi_1(x)\in C_r,$$
{where $r\in \{\{1,\cdots,R\}: \Phi_1(\bar x)\in C_r\}$.}
  Since the abstract set in the subsystem is a convex polyhedral set, we can propose suitable conditions on problem data such that RCPLD is a sufficient condition for the existence of a local  error bound for each subsystem. Since the number of subsystem is finite, the local error bound for the original system also exists provided that the one for each subsystem exists.
%
 
%

We summarize the main contributions as follows:
\begin{itemize}
\item Using the generator approach, we propose a novel concept of RCPLD for the disjunctive system which is weaker than the   one  proposed  in Xu and Ye \cite{xy2}. 
We apply the new RCPLD to the ortho-disjunctive program to obtain RCPLD for the ortho-disjunctive program. Such a condition  is  a unifying framework for  recovering  the existing RCPLD for certain classes of mathematical programs, such as  MPEC-RCPLD and MPSC-RCPLD.
 It also provides a  tool to obtain a new RCPLD condition for  disjunctive systems such  as MPVC.

\item 
By using the local decomposition approach, we propose the concept of  piecewise RCPLD. We show that when each disjunctive set $\Gamma_i$ is a subset of the space $\mathbb{R}^{p_i}$ with $p_i=1, 2$, the piecewise RCPLD implies RCPLD. Moreover we give an example to show that this implication fails to hold if $p_i\geq 3$ for at least one $i$. The piecewise RCPLD ensures the existence of a local error bound for any system from  applications that can be reformulated as a disjunctive system; for example, KKT system for a nonlinear program ( \cite{f02,ms,r82}). We also observe that if RCPLD coincides with the piecewise RCPLD, then RCPLD leads to error bounds as well. This is true for MPSC. We observe that MPSC-RCPLD coincides with the MPSC piecewise RCPLD and hence 
it is a constraint qualification, which also leads to the error bound property and the exact penalization.

\item We introduce   MPDC variants of certain constraint qualifications that are stronger than the RCPLD such as the constant rank constraint qualification (CRCQ), the relaxed constant rank constraint qualification (RCRCQ), the constant positive linear dependence constraint qualification (CPLD), the enhanced relaxed constant positive linear dependence condition (ERCPLD) and study their relationships with RCPLD and  the piecewise RCPLD.  
Applying them to the ortho-disjunctive systems, we {can} derive weaker MPVC variants of CRCQ, CPLD and MPSC variants of CRCQ, RCRCQ than those defined in the literature (\cite{hks12} for MPVC and \cite{lg} for MPSC),  and new MPVC variants of RCRCQ,  and ERCPLD, which were not introduced before.


\end{itemize}

The paper  is organized as follows.  In Section 2, we provide  the preliminaries including the representation for the normal cones of the disjunctive set. In Section 3, we introduce RCPLD for the disjunctive system and show that it is a constraint qualification  and derives M-stationarity condition.  
{In Section 4, by using the local decomposition approach, the piecewise RCPLD is proposed to obtain the error bound property for MPDC.
Moreover, we review and introduce some MPDC variants of  constraint qualifications    and study their relationships with RCPLD and piecewise RCPLD.}
In Section 5, we apply  RCPLD and the piecewise RCPLD to the ortho-disjunctive programs, including MPEC, MPVC and MPSC. The 
concluding remarks can be found in the final section.


We now explain our notation.
For $g: \mathbb{R}^d\rightarrow \mathbb{R}$,  the
 gradient vector at $x$ is a column vector denoted by $\nabla g(x)$ and the maximum part of $g$ at $x$ is denoted by $g^+(x):=\max\{0, g(x)\}$.
For sets $C_1, \cdots, C_l$, we denote by $\Pi_{i=1}^I C_i$ their Cartesian {product}.
A vector $e_i\in \mathbb{R}^m$ is the unit vector where $i$-th component is one.
{Denote by $\|\cdot\|$ and $\|\cdot\|_2$ any  norm and 2-norm in $\mathbb{R}^n$, respectively. }
For a set of  finitely many vectors $A:=\{ v_1,\cdots, v_n\}$, we denote by
${\rm span} (A)$ and ${\rm cone} (A)$ the subspace and the convex cone generated by the vectors in $A$, respectively.
Given   finite index sets $I, J$,  a pair  $(\{v_i\}_{i\in I}, \{u_i\}_{i\in J})$ of family of vectors $\{v_i\}_{i\in I}, \{u_i\}_{i\in J}\subseteq \mathbb{R}^d$ is said to be positive linearly dependent if there exist scalars $\{\alpha_i\}_{i\in I}$ and $\{\beta_i\}_{i\in J}$ with $\alpha_i\geq 0$ for any $i\in I$, not all equal to zero such that $\sum_{i\in I}\alpha_i v_i+\sum_{i\in J}\beta_i u_i=0$. 

\section{Preliminaries}
We first give notations to subdifferentials and normal cones.  
We omit their definitions and refer the readers to the standard reference of variational analysis
 in \cite{c,clsw,m1,var}.
For  a Lipschitz continuous function $\phi: \mathbb{R}^d \rightarrow  \mathbb{R}$ at $\bar  x$,
 denote by $\widehat{\partial} \phi( \bar x)$ and $\partial \phi(\bar{x})$
 the Fr\'{e}chet/regular  subdifferential and  the limiting/Mordukhovich/basic  subdifferential
 of $\phi $ at $ \bar  x$, respectively.
$\phi$ is subdifferentially regular at $\bar x$ if $\partial \phi(\bar{x})=\widehat{\partial} \phi( \bar  x)$ (\cite[Corollary 8.11, Theorem 9.13]{var}).
For a closed set $C\subseteq \mathbb{R}^d$, denote the Fr\'{e}chet/regular normal cone and the limiting/Mordukhovich/basic normal cone at $\bar x\in C$    by $\widehat{\mathcal{N}}_{C}(\bar x)$
and $\mathcal{N}_{C}(\bar x)$, respectively.
When $C$ is convex, all the normal cones are equal to the normal cone in the sense of convex analysis. A closed set $C$ is subdifferentially/Clarke regular at $\bar x$ if  $\mathcal{N}_{C}(\bar x)=\widehat{\mathcal{N}}_{C}(\bar x)$ (\cite[Definition 6.4]{var}). 

The following  result which is an extension of  Carath\'{e}odory's lemma \cite[Theorem 17.1]{con} will be useful.
\begin{lemma}\cite[Lemma 1]{RCPLD}\label{lem3-1}
If $v=\displaystyle \sum_{i=1}^{m+n} \alpha_i v_i$ with $v_i\in\mathbb{R}^d$ for every $i$, $\{v_i\}_{i=1}^m$ is linearly independent and $\alpha_i\neq 0$, $i=m+1,\cdots,m+n$, then there exist $I\subseteq \{m+1,\cdots,m+n\}$ and scalars $\bar{\alpha}_i$ for every $i\in \{1,\cdots,m\}\cup I$ such that\\
(i) $v=\displaystyle \sum_{\{1,\cdots,m\}\cup I} \bar{\alpha}_i v_i$ with $\alpha_i \bar{\alpha}_i>0$ for every $i\in I$,\\
(ii) $\{v_i\}_{\{1,\cdots,m\}\cup I}$ is linearly independent.
\end{lemma}
It is clear that  the vector $v$ in  Lemma \ref{lem3-1}  can not be a zero vector.

For any two subsets $A^I,A^E$ in $\mathbb{R}^d$, denote by
\begin{eqnarray*}
{\cal G}(A^I, A^E):={\rm cone}(A^I)+{\rm span}(A^E).
\end{eqnarray*}
Given a convex polyhedral cone $D\subseteq \mathbb{R}^d$.  A pair of family of vectors $(A^I,A^E)$ is said to be a {\it generator} of $D$,   {if the sets $A^I,A^E$ are finite and disjoint, the vectors in $A^E$ are linearly independent, 
$A^I\cup A^E$ has the minimal number of vectors such that}
\begin{equation}\label{genercone} 
D={\cal G}(A^I, A^E).
\end{equation}
{Note that if $D=\{0\}$, then $A^I\cup A^E=\{0\}$.}
We use the following example to explain the concept of a generator for a convex polyhedral cone.
\begin{example} Consider {$D=\mathbb{R}^2$}.  
Then $(A^I,A^E)$ with $A^I:=\emptyset, A^E:=\{(1,0),(0,1)\}$ is a generator for $D$. But $(\widetilde{A}^I,\widetilde{A}^E)$ with $\widetilde{A}^I:=\{(1,1),(-1,1)\}, {\widetilde{A}^E:=\{(0,1)\}}$ is not a generator for $D$ since although (\ref{genercone}) holds and the vectors in $A^E$ are linearly independent, the set $\widetilde{A}^I\cup \widetilde{A}^E=\{(1,1),(-1,1),(0,1)\}$ is not minimal since it has one more  vector than the previous one. 
\end{example} 

%

We now describe the representation for the normal cones of the disjunctive set in the form of   $\Gamma=\bigcup_{r=1}^R C_r\subseteq \mathbb{R}^d$ with $C_r$ being a convex polyhedral set, $r=1,\cdots, R$. Suppose each $C_r$ has the following representation
\begin{eqnarray*}
C_r:=\left \{x\in \mathbb{R}^d:  \langle c_j^r, x\rangle\leq \alpha_j^r,\ j\in {\cal I}_r , \langle c_j^r, x\rangle= \alpha_j^r,   j\in {\cal E}_r\right \},
\end{eqnarray*}   
where ${\cal I}_r , {\cal E}_r$ are finite index sets,  $0\not =c_j^r\in \mathbb{R}^d$ and $\alpha_j^r\in \mathbb{R}$ are given constants. Also we choose the generator for equality constraints $A_r^E:=\{c_j^r: j\in {\cal E}'_r\}$, {where ${\cal E}'_r\subseteq {\cal E}_r$,} to be linearly independent.
For $x \in C_r$, let ${\cal I}_r(x):=  \{j\in {\cal I}_r: \langle c_j^r, x\rangle= \alpha_j^r  \}$ be the index set of active constraints and  $A_r^I(x):=\{c_j^r: j\in {\cal I}_r(x)\}$.
Without loss of generality, we assume that for any $j,j'\in {\cal I}_r(x)$, $c_j^r\neq -c_{j'}^r$, otherwise we can move $c_j^r$ to the set $A_r^E$. 
{By  \cite[Theorem 6.46]{var}, for each $r$ and any $x\in C_r$, the normal cone is convex polyheral and has the representation:
$${\cal N}_{C_r} (x)=\left \{ \sum_{j\in {\cal I}_r(x)} \lambda_j c_j^r+ \sum_{j\in {\cal E}_r} \lambda_j c_j^r: \lambda_j \geq 0, j\in {\cal I}_r(x)\right \}={\cal G}(A_r^I(x), A_r^E).$$ Since the normal cone of $C_r$ is a convex polyheral cone, we can find  a generator for it.
  Although  this generator may be different, for simplicity we do not change the notation and still denote a generator of ${\cal N}_{C_r} (x)$ by
$(A_r^I(x), A_r^E)$. We also denote the set
 $${A}_{C_r}(x):=A_r^I(x)\cup A_r^E.$$

For a point $x\in \Gamma=\bigcup_{r=1}^R C_r$,  denote the active index set for $\Gamma$  at $x$ as $I(x):=\{r: x\in C_r\}. $ Then for $x\in \Gamma$, from Mehlitz \cite[Lemma 2.2]{mmpscnew},
\begin{equation}\label{regularcone} 
\widehat{\mathcal{N}}_{\Gamma}(x)=\bigcap_{r\in I(x)}{\mathcal{N}}_{C_r}(x).
\end{equation}
Since $\widehat{\mathcal{N}}_{\Gamma}(x)$ is a convex polyhedral cone, we can find a generator for it. We denote a generator of $\widehat{\mathcal{N}}_{\Gamma}(x)$  by  $(\widehat{A}_{\Gamma}^I(x), 
\widehat{A}_{\Gamma}^E(x))$ and the set 
$$\widehat{A}_{\Gamma}(x):=\widehat{A}_{\Gamma}^I(x)\cup \widehat{A}_{\Gamma}^E(x) .$$
 Hence we have
\begin{equation}
 \widehat{\mathcal{N}}_{\Gamma}( x)={\cal G} (\widehat{A}_{\Gamma}^I(x), \widehat{A}_{\Gamma}^E(x)).
 \label{Ahat}
\end{equation}

From (13) and (14) in Adam et al. \cite{acp},  for each $\bar x\in \Gamma$, we can find $\delta>0$ such that
\begin{equation}\label{limitingcone} 
\mathcal{N}_{\Gamma}(\bar x)=\bigcup_{x\in \mathbb{B}_{\delta}(\bar x)} \widehat{\mathcal{N}}_{\Gamma}(x).
\end{equation}
Define the sets
\begin{equation}
A_\Gamma^I(\bar x):=\bigcup_{x\in \mathbb{B}_{\delta}( \bar x)} \widehat{A}_{\Gamma}^I(x), \qquad A_\Gamma^E(\bar x):=\bigcup_{x\in \mathbb{B}_{\delta}( \bar x)} \widehat{A}_{\Gamma}^E(x), \quad A_\Gamma(\bar x):=A_\Gamma^I(\bar x)\cup A_\Gamma^E(\bar x),
\label{Ax}
\end{equation}
 where $\delta$ is the constant satisfying condition {(\ref{limitingcone})}.  We call 
{$(A_\Gamma^I(\bar x), A_\Gamma^E(\bar x))$}
  {\it the generator set for the limiting normal cone of $\Gamma$ at $\bar x$}.
 Note that since  $\mathcal{N}_{\Gamma}(\bar x)$ is not  convex  in general, 
a generator set for the  limiting normal cone is not a generator in the sense of (\ref{genercone}).
 
 \begin{remark}\label{nor-lim}
By (\ref{limitingcone}),  there exists a sufficiently small {$\delta'<\delta$} such that for any $x\in \mathbb{B}_{\delta'}(\bar x)\cap \Gamma$,
\begin{eqnarray}
\mathcal{N}_{\Gamma}(x)\subseteq \mathcal{N}_{\Gamma}(\bar x).\label{normalcone}
\end{eqnarray}
Indeed,  by (\ref{limitingcone}) for each $x\in {\mathbb{B}_{\delta'}(\bar x)\cap }\Gamma$ there exists $\delta_x>0$ sufficiently small such that $\mathbb{B}_{\delta_x}(x)\subseteq  \mathbb{B}_{\delta}(\bar x)$ 
and $$
\mathcal{N}_{\Gamma}(x)=\bigcup_{y\in \mathbb{B}_{\delta_x}(x)} \widehat{\mathcal{N}}_{\Gamma}(y)\subseteq \bigcup_{y\in \mathbb{B}_{\delta}(\bar x)} \widehat{\mathcal{N}}_{\Gamma}(y)=\mathcal{N}_{\Gamma}(\bar x)\quad  \mbox{ for all } x\in \mathbb{B}_{\delta'}(\bar x),
$$ and hence (\ref{normalcone}) holds.
Similarly, from (\ref{Ax}) and (\ref{Ahat}) we have 
\begin{equation}
 \widehat{A}_{\Gamma}^I(y) \subseteq A_{\Gamma}^I(\bar x), \quad \widehat{A}_{\Gamma}^E(y) \subseteq A_{\Gamma}^E(\bar x),
 \quad {\cal G} (\widehat{A}_{\Gamma}^I(y), \widehat{A}_{\Gamma}^E(y))\subseteq \mathcal{N}_{\Gamma}(\bar x),\quad 
{\forall y\in \mathbb{B}_{\delta_x}(x)} \label{eqn2.8}
\end{equation}
 and
\begin{eqnarray}
A_{\Gamma}^I(x) \subseteq A_{\Gamma}^I(\bar x), \quad A_{\Gamma}^E(x) \subseteq A_{\Gamma}^E(\bar x)\quad \forall x\in \mathbb{B}_{\delta'}(\bar x).\label{Ainclusion}
\end{eqnarray} 
\end{remark}

From Mehlitz \cite[Lemma 2.2]{mmpscnew},  we have
\begin{eqnarray}\label{unionnorm}
\displaystyle\mathcal{N}_\Gamma(\bar x)\subseteq \displaystyle\bigcup_{r\in I(\bar x)}  {\mathcal{N}}_{C_r} (\bar x).
\end{eqnarray}
Hence it is interesting to know whether or not  the elements in the generator set $A_\Gamma(\bar x)$ must be selected from the union of all $ A_{C_r}(\bar x)$ for $r\in I(\bar x)$.
 In the following lemma, we show that the answer is affirmative if the dimension $d$ is less or equal to two.   It will be shown in Example \ref{Counterexample} that the answer is negative if the dimension $d$  is greater than two. The reason is that (\ref{gernator-1}) may not hold if $d$ is greater than two.

\begin{lemma}\label{Lem2.1} Let $\Gamma=\bigcup_{r=1}^R C_r \subseteq \mathbb{R}^d$ with $d\leq2$ and $C_r$ being convex polyhedral.  Then for any $\bar x\in \Gamma$, we have
$$A_\Gamma (\bar x)\subseteq \bigcup_{r\in I(\bar x)}  A_{C_r}(\bar x).$$
\end{lemma}
{\bf Proof.} 
By (\ref{Ainclusion}),  for any $x\in \mathbb{B}_{\delta}(\bar x)$ we have
 \begin{equation}\label{gernator-1h}
    A_{C_r}(x)\subseteq A_{C_r}(\bar x).
 \end{equation} 

Since $d\leq 2$ and by (\ref{regularcone}), 
$\widehat{\mathcal{N}}_{\Gamma}(x)$ is the intersection of finitely many convex cones ${\mathcal{N}}_{C_r}(x)$,  we have
 \begin{equation} 
 \widehat{A}_\Gamma(x) \subseteq \bigcup_{r\in I(x)} 
   A_{C_r}(x).\label{gernator-1}
  \end{equation}
Furthermore, 
 from (\ref{Ax}) and (\ref{gernator-1h})-(\ref{gernator-1}),
\begin{eqnarray*}
A_\Gamma(\bar x)&=&\bigcup_{x\in \mathbb{B}_{\delta}(\bar x)} \widehat{A}_\Gamma(x)\subseteq \bigcup_{x\in \mathbb{B}_{\delta}(\bar x)} \bigcup_{r\in I(x)}   A_{C_r}(x)\subseteq \bigcup_{r\in I(\bar x)}   A_{C_r}(\bar x),
\end{eqnarray*}
where the last inclusion followed by the closedness of $C_r$.
\BOX

We now review the concept of M-stationarity below.  The reader is referred to \cite{bg,fko} for other notions of stationarity for MPDC such as Strong-stationarity, Q-stationarity, etc.
Recall that ${\cal F}$ is the solution set of the disjunctive system $(\ref{feasibility})$.
\begin{defn}[Mordukhovich {stationarity}  condition](\cite[Definition 4]{bg}
)\label{Mstationarypoint}
Consider  a feasible solution $\bar x$ of the disjunctive program:
$ (P) \ \min f(x) \ \mbox{ s.t. } x\in {\cal F},$
where $f$ is Lipschitz continuous at $\bar x$.
We call $\bar x$  an M-stationarity point  if there exists    $(\lambda^g,\lambda^h, \bar \eta) \in \mathbb{R}^{n}\times\mathbb{R}^{m}\times \Pi_{i=1}^l \mathbb{R}^{ p_i}$  such that 
 \begin{eqnarray*}
&& 0\in \partial f(\bar x)+\sum_{i\in \bar I_g} \lambda_i^g \partial g_i(\bar x) +\sum_{i=1}^m \lambda_i^h \nabla h_i(\bar x) +\sum_{i=1}^l \nabla \Phi_i(\bar x)^T \bar \eta_i,\\
 &&\lambda_i^g\geq 0, i\in \bar I_g,\ \bar \eta_i\in \mathcal{N}_{\Gamma_i}(\Phi_i(\bar x)), \ \  i=1,\cdots, l,
\end{eqnarray*}
where  $\bar{I}_g:=I_g(\bar x):=\{i=1,\cdots,n: g_i(\bar x)=0\}$. 
\end{defn}


\begin{defn}[Error bound property]\label{Defnerr}
We say that  $\bar   x \in {\cal F}$  satisfies the {\rm error bound property}  for system (\ref{feasibility}) if we can find $\alpha\geq 0$ and $\varepsilon>0$ such that 
$$
d_{\cal F}(x) \leq \alpha \left ( \|g^+(x)\| + \|h(x)\| + \displaystyle \sum_{i=1}^l d_{\Gamma_i}(\Phi_i(x))\right ), \quad \forall x\in \mathbb{B}_{\varepsilon}(\bar x),
$$
where 
$d_{\cal F}(x)$ is the distance  from $x$ to $\cal F$.
\end{defn}

\section{RCPLD for disjunctive systems}
We propose the following RCPLD condition for the disjunctive system (\ref{feasibility}) and show that it is a constraint qualification for {M-stationarity. }
\begin{defn}[{MPDC-RCPLD}]\label{nrcpld} 

We say $\bar x\in {\cal F}$ satisfies the MPDC-RCPLD for the disjunctive system (\ref{feasibility}) if conditions (i)-(ii) are satisfied.
\begin{itemize} 
\item[{\rm (i)}] $\{\nabla h_i(x)\}_{i=1}^m$ has constant rank,  for any $x\in U(\bar x)$,  a neighborhood of $\bar x$. 
\item[{\rm (ii)}] Let $J\subseteq \{1,\cdots,m\}$ be the index set that $\{\nabla h_i(\bar x)\}_{i\in J}$  is a basis of ${\rm span} \{\nabla h_i(\bar x)\}_{i=1}^m.$
Suppose   there exist $I\subseteq \bar I_g$, 
 a nonzero vector $(\lambda^h,\lambda^g,\bar \eta) \in \mathbb{R}^{m}\times \mathbb{R}^{n}\times \Pi_{i=1}^l \mathbb{R}^{ p_i}$ with $\lambda_i^g\geq 0$,  $v_i \in \partial g_i(\bar x), i\in I$ and $\bar \eta_i\in \mathcal{N}_{{\Gamma_i}}(\Phi_i(\bar x)),  i\in \{1,\cdots,l \}$ such that
\begin{eqnarray}
&& 0 = \sum_{i\in I} \lambda_i^g v_i+\sum_{i\in J} \lambda_i^h \nabla h_i(\bar x)  + \sum_{i=1}^l \nabla \Phi_i(\bar x)^T \bar  \eta_i.\label{nna1}
\end{eqnarray}  

Then for sufficiently large $k$,
\begin{eqnarray*}
\{v_i^k\}_{i\in I}\cup\{\nabla h_i(x^k)\}_{i\in J}\cup  
\bigcup_{\beta_i\in A_i^I\cup A_i^E,i\in \{1,\cdots,l \}} \{\nabla \Phi_i(x^k)^T \beta_i\} \label{normalvectorphi}
\end{eqnarray*} 
is linearly dependent, for any
$\{x^k\}, \{v^k\}$ such that $x^k\rightarrow \bar x$, $x^k\neq \bar x$,  
 $v^k_i\in \partial g_i(x^k)$, $v^k_i\rightarrow v_i$,
any set of 
 linearly independent 
vectors $A_i:=A_i^I\cup A_i^E$ where $A_i^I \subseteq A_{\Gamma_i}^I(\Phi_i(\bar x))$, $A_i^E \subseteq A_{\Gamma_i}^E(\Phi_i(\bar x))$ satisfying  $0\not =\bar \eta_i\in {\cal G}{(A_i^I,A_i^E)} \subseteq \mathcal{N}_{{\Gamma_i}}(\Phi_i(\bar x))$ and $A_i:=A_i^I\cup A_i^E=\emptyset$ if $\bar \eta_i=0$.
\end{itemize}
\end{defn}

We are now ready to show that MPDC-RCPLD  is a constraint qualification.
\begin{thm}\label{OC2}
Suppose $\bar x$ is a local solution of 
$ ({\rm P}): \min_{x} f(x) \ \mbox{ s.t. } x\in {\cal F},$ where $f$ is Lipschitz continuous at $\bar x$.
 Suppose that {MPDC-RCPLD}  holds at $\bar x$. 
Then $\bar x$ is a Mordukhovich stationary point for  $(\rm P)$ as defined in Definition \ref{Mstationarypoint}.
\end{thm}
{\bf Proof.} 
To concentrate on the main idea,  in the proof we assume that $l=1$, $\Gamma_1=\Gamma \subseteq \mathbb{R}^p$ and $\Phi(x)=\Phi_1(x)$.  We also omit the equality and inequality constraints and refer the reader to
\cite[Theorem 3.1]{xy2} for the corresponding treatments.   

We introduce an auxiliary variable $y$ and consider the problem 
\begin{eqnarray*}
({\rm \widetilde P})~~~~~~~\min_{x,y}&&  f(x)\nonumber\\
{\rm s.t.} && 
\Psi(x,y):=\Phi(x)-y=0,  (x,y) \in  \mathbb{R}^d\times {\Gamma}.
\label{afeasibility}
\end{eqnarray*}
 If $\bar x$ is a local solution for (P), then  $(\bar x,\bar y)$ with $\bar y=\Phi(\bar x)$ is a local solution of $({\rm  \widetilde{P}})$ and we can find $\varepsilon>0$ such that for any $(x,y)\in  {\bar \mathbb{B}}_{\varepsilon}(\bar{x},\bar y)$, $(x,y)$ is feasible for $({\rm  \widetilde{P}})$, we have that $f(\bar{x})\leq f(x)$, where ${\bar \mathbb{B}}_{\varepsilon}(\bar{x},\bar y)$  is the closure of the open  ball {centered} at $(\bar{x},\bar y)$ with  radius equal to $\varepsilon$.

Step 1:  
Consider the penalized problem:
\begin{eqnarray*}
(\widetilde{\rm P}_k)~~~\min && G_k(x,y):=f(x)+\frac{k}{2}\sum_{i=1}^ p \Psi_i^2(x,y)+\frac{1}{2}\|(x,y)-(\bar{x},\bar y)\|^2_2 \\
 {\rm s.t.} && (x,y)\in (\mathbb{R}^d\times {\Gamma})\cap { \bar \mathbb{B}}_{\varepsilon}(\bar{x},\bar y).
\end{eqnarray*}

For any $k$, there exists an solution $(x^k,y^k)$ of $(\widetilde{\rm P}_k)$. 
Without loss of generality, assume $\displaystyle \lim_{k\to\infty}(x^k,y^k)=(x^*,y^*)$, then $(x^*,y^*)\in (\mathbb{R}^d\times {\Gamma})\cap { \bar \mathbb{B}}_{\varepsilon}(\bar{x},\bar y)$.
Since $(x^k,y^k)$ is the minimizer of $(\widetilde{\rm P}_k)$, 
\begin{eqnarray}\label{th2-2-1}
f(x^k)+\frac{k}{2}\sum_{i=1}^ p \Psi_i^2(x^k,y^k)+\frac{1}{2}\|(x^k,y^k)-(\bar{x},\bar y)\|^2_2
= G_k(x^k,y^k) \leq G_k(\bar{x},\bar y)=f(\bar{x}).
\end{eqnarray}
From the continuity of $f$, the sequence $\left \{f(x^k) \right \}$ is bounded. Thus we have that
$
\Psi(x^k,y^k)\to 0=\Psi(x^*,y^*)
$
as $k\to\infty$. 
Therefore $(x^*,y^*)$ is feasible for $({\rm  \widetilde{P}})$.  

Condition $(\ref{th2-2-1})$ implies that
$
f(x^k)+\frac{1}{2}\|(x^k,y^k)-(\bar{x},\bar y)\|^2_2 \leq f(\bar{x}).
$
Taking limit as $k\to\infty$,  
$
f(x^*)+\frac{1}{2}\|(x^*,y^*)-(\bar x,\bar y)\|^2_2 \leq f(\bar{x})\leq f(x^*),
$
then $(x^*,y^*)=(\bar x,\bar y)$ and  $\{(x^k,y^k)\} $ converges to $(\bar x,\bar y)$. 

Step 2: It is obvious that $(x^k,y^k)$ is an interior point of ${\bar \mathbb{B}}_{\varepsilon}(\bar{x},\bar y)$ for large $k$. From the necessary optimality condition
 and the nonsmooth sum and chain rules, 
(see e.g. Rockafellar and Wets \cite[Theorem 8.15, Theorem 10.6, Exercise 10.10]{var} and Mordukhovich \cite[Theorem 2.33]{m1}), 
there exist
 $v_0^k\in \partial f(x^k)$, 
$u^k\in k \sum\limits_{i=1}^p \Psi_i(x^k,y^k) \nabla \Psi_i(x^k,y^k)+\{0\}^d\times \mathcal{N}_{\Gamma}(y^k)$ such that
\begin{eqnarray}
0=v_0^k\times\{0\}^p+u^k+ \left ((x^k,y^k)-(\bar{x},\bar y)\right ).\label{sumrule}
\end{eqnarray}
For large $k$, we assume that $u^k\neq 0$, since otherwise by taking limit {in} $(\ref{sumrule})$ as $k\to \infty$, $0\in \partial f(\bar x)$ 
from the outer semi-continuity \cite[Proposition 8.7]{var} and then $\bar x$ is an M-stationary point automatically. 

Denote by $\Phi^i(x)$   the $i$th component of $\Phi(x)$. Then
 there exist 
$\lambda_{i,k}^{\Phi},i=1,\cdots,p$ and $\eta^k\in \mathcal{N}_{\Gamma}(y^k)$ such that
\begin{eqnarray*}
u^k= 
 \sum_{i=1}^p \lambda_{i,k}^{\Phi} \left (\nabla \Phi^i(x^k)\times (-e_i) \right )
 +\{0\}^d\times \eta^k.
\end{eqnarray*}

From $(\ref{limitingcone})$, 
similarly as in Remark \ref{nor-lim}, 
there exist $\delta>0$ and small enough $\delta^k>0$
such that 
$\mathbb{B}_{\delta^k}(y^k)\subseteq \mathbb{B}_{\delta} (\bar y)$,
$$
\mathcal{N}_{\Gamma}(y^k)=
\cup_{y\in \mathbb{B}_{\delta^k} (y^k)} \widehat{\mathcal{N}}_{\Gamma}(y)=
\cup_{y\in \mathbb{B}_{\delta^k}(y^k)} {\cal G}(\widehat{A}_\Gamma^I(y), \widehat{A}_\Gamma^E(y))
\subseteq \mathcal{N}_{\Gamma}(\bar y).
$$ 
Hence for 
$\eta^k\in \mathcal{N}_{\Gamma}(y^k)$, there exists $\tilde{y}^k\in \mathbb{B}_{\delta^k} (y^k)$ such that 
$$\eta^k\in \widehat{\mathcal{N}}_{\Gamma}(\tilde{y}^k)={\cal G}(\widehat{A}_\Gamma^I(\tilde{y}^k), \widehat{A}_\Gamma^E(\tilde{y}^k)).$$
Denote by 
$(\widehat{A}_{k}^I, \widehat{A}_{k}^E):=(\widehat{A}_\Gamma^I(\tilde{y}^k), \widehat{A}_\Gamma^E(\tilde{y}^k))$. Then we have 
\begin{equation} 
u^k\in  \sum_{i=1}^p \lambda_{i,k}^{\Phi} \left (\nabla \Phi^i(x^k)\times (-e_i) \right )
+\{0\}^d\times {\cal G} (\widehat{A}_{k}^I, \widehat{A}_{k}^E)\label{calculus}
\end{equation} 
and by (\ref{eqn2.8}) we have
\begin{equation}
 \widehat{A}_{k}^I \subseteq A_{\Gamma}^I(\bar y), \quad \widehat{A}_{k}^E \subseteq A_{\Gamma}^E(\bar y),
 \quad {\cal G} (\widehat{A}_{k}^I, \widehat{A}_{k}^E)\subseteq \mathcal{N}_{\Gamma}(\bar y).\label{generatork}
\end{equation} 
It is easy to see that for any $k$, 
the set of vectors $\large \{\nabla \Phi^i(x^k)\times (-e_i)\large \}_{i=1}^p$ is linearly independent.
 By Lemma \ref{lem3-1}, we obtain subsets 
$\widetilde{A}_{k}^I\subseteq\widehat{A}_{k}^I, \widetilde{A}_{k}^E\subseteq\widehat{A}_{k}^E$
and $\bar{\lambda}_{k}^{\Phi}$
such that
\begin{eqnarray}\label{sum1}
&&u^k\in 
\sum_{i=1}^p \bar{\lambda}_{i,k}^{\Phi} \left (\nabla \Phi^i(x^k)\times (-e_i) \right )+\{0\}^d\times {\cal G} (\widetilde{A}_{k}^I, \widetilde{A}_{k}^E),\\
&& {\cal G} (\widetilde{A}_{k}^I, \widetilde{A}_{k}^E)\subseteq {\cal G} (\widehat{A}_{k}^I, \widehat{A}_{k}^E)
\label{Obviously} 
\end{eqnarray}
and 
\begin{eqnarray}\label{lin}
\left \{\nabla \Phi^i(x^k)\times (-e_i) \right \}_{i=1}^p
\cup \{\{0\}^d\times \beta\}_{\beta\in \widetilde{A}_{k}^I \cup\widetilde{A}_{k}^E}
\end{eqnarray}
 is linearly independent. 

 Since the sets $\widehat{A}_{k}^I$ and $\widehat{A}_{k}^E$ are both finite for every large $k$, taking a subsequence if necessary, we  may assume that
 $\widetilde{A}_{k}^I\equiv A^I$, $\widetilde{A}_{k}^E\equiv A^E$. 
 {Then from (\ref{generatork}) and (\ref{Obviously}), we have that ${\cal G} (A^I, A^E) \subseteq \mathcal{N}_{\Gamma}(\bar y)$ and  $A^I \subseteq A_{\Gamma}^I(\Phi(\bar x))$, $A^E \subseteq A_{\Gamma}^E(\Phi(\bar x))$. }
 
 From the conditions (\ref{sumrule}) and (\ref{sum1}), we obtain that 
 \begin{eqnarray}\label{sumrule1new}
&&0 \in v_0^k\times\{0\}^p+(x^k,y^k)-(\bar{x},\bar y)
+ \sum_{i=1}^p \bar{\lambda}_{i,k}^{\Phi} \left  (\nabla \Phi^i(x^k)\times (-e_i)\right )+\{0\}^d\times {\cal G} (A^I, A^E),\\
&&~~~
\{\nabla \Phi(x^k)^T \beta\}_{\beta\in A^I \cup A^E} \mbox{ is linearly independent,} \label{addcondition}
\end{eqnarray}
where the condition (\ref{addcondition}) follows from the linear independence of the set of vectors in (\ref{lin}).

Note that if $x^k=\bar x$, then $(\ref{th2-2-1})$ implies that $\frac{1}{2}\|(\bar x,y^k)-(\bar{x},\bar y)\|_2^2\leq 0$ and thus $(x^k,y^k)=(\bar{x},\bar y)$. Then condition \((\ref{sumrule1new})\) implies the M-stationary condition. Thus  we assume that $x^k\neq\bar x$.

Step 3: If the sequence $\{\bar{\lambda}_{k}^{\Phi}\}$
is bounded, then we assume $\bar{\lambda}_{k}^{\Phi}\to \bar{\lambda}^{\Phi}$ as $k\rightarrow \infty$ without loss of generality.
Since $f$ is  locally Lipschitz continuous at $\bar x$ and $x^k \to x$, we have that $v_0^k$ is bounded and the limits of $\{v_0^k\}$ belongs to $\partial f(\bar x)$
 by the outer semi-continuity of the limiting subdifferential. 
Taking limits in \((\ref{sumrule1new})\), we have
\begin{eqnarray*}
 0\in \partial f(\bar x)
 \times\{0\}^p  + 
 \sum_{i=1}^p \bar{\lambda}_{i}^{\Phi} \left (\nabla \Phi^i(\bar x)\times (-e_i)\right )+\{0\}^d\times\mathcal{N}_{\Gamma}(\Phi(\bar x)),
\end{eqnarray*}
from which we obtain the stationary condition and the proof is therefore completed.

Otherwise if $\{\bar{\lambda}_{k}^{\Phi}\}$ is not bounded 
and there exists a subsequence $K$, 
$\|\bar{\lambda}_{k}^{\Phi}\|_2\to\infty$ as $ k\rightarrow \infty, k \in K$, then without loss of generality, we assume that 
$$\displaystyle\lim_{k\to\infty,k\in K}\frac{\bar{\lambda}_{k}^{\Phi}}{\|\bar{\lambda}_{k}^{\Phi}\|_2}=\lambda^{\Phi}.$$ 
Assume 
$v_0^k\to v_0\in \partial f(\bar x)$ by the outer semi-continuity of the limiting subdifferential. 
Dividing by $\|\bar{\lambda}_{k}^{\Phi}\|_2$ on \((\ref{sumrule1new})\) and letting $k \rightarrow \infty, k\in K$,
we obtain that
 $0\neq\bar \eta\in {\cal G} (A^I, A^E)\subseteq \mathcal{N}_{\Gamma}(\bar y)$ with $A^I\cup A^E$ being linearly independent from (\ref{addcondition}) such that
\begin{eqnarray*}
0 =
\sum_{i=1}^p \lambda_i^{\Phi} \left (\nabla \Phi^i(\bar x)\times (-e_i) \right )
+\{0\}^d\times\bar \eta,
\end{eqnarray*}
it implies that
$0 = 
\nabla \Phi(\bar x)^T \bar{\eta}.
$
From  MPDC-RCPLD  (ii),  $
\{\nabla \Phi(x^k)^T \beta\}_{\beta\in A^I \cup A^E}$
is  linearly dependent if $\bar \eta\neq 0$.
 But this  contradicts (\ref{addcondition}). 
The contradiction proves that  $\{\bar{\lambda}_{k}^{\Phi}\}$ is bounded. 
\BOX

Similar to the proof of Xu and Ye \cite[Proposition 3.1]{xy2}, we can prove that {MPDC-RCPLD} is persistent locally.
\begin{proposition}\label{persist}  
{Assume that  {MPDC-RCPLD} for the disjunctive system (\ref{feasibility}) holds at $\bar x\in \cal{F}$. Then {MPDC-RCPLD}  satisfies at any  point in a small neighborhood  of $\bar x$. }
\end{proposition}   
{\bf Proof.}  Consider any  sequence $x^k\rightarrow \bar x$ as $k\to\infty$.  
It is obvious that MPDC-RCPLD (i) holds at each $x^k$ when $k$ is sufficiently large.
We only need to prove the condition (ii). 
Suppose $J\subseteq \{1,\cdots,m\}$ is the index set such that $\{\nabla h_i(\bar x)\}_{i\in J}$  is a basis of ${\rm span} \{\nabla h_i(\bar x)\}_{i=1}^m.$
Then $\{\nabla h_i(x^k)\}_{i\in J}$ is linearly independent and thus a basis of span$\{\nabla h_i(x^k)\}_{i=1}^m$ for any large $k$ from MPDC-RCPLD (i).

To the contrary, assume that  MPDC-RCPLD(ii) does not hold at each point of a subsequence $\{x^k\}_{k\in K}$. 
Then there exist $\mathcal{I}^k\subseteq I_g(x^k)$, a nonzero vector $(\lambda^g_{k},\lambda^h_{k}, \eta^k )$ with $\lambda^g_{i,k}\geq 0$, $v_i^k\in\partial g_i(x^k)$,  $i\in \mathcal{I}^k$ and $\eta^k_i\in\mathcal{N}_{{\Gamma_i}}(\Phi_i(x^k))$, $i=1,\cdots, l$ such that 
\begin{eqnarray}\label{persisnew}
0 = \sum_{i\in \mathcal{I}^k} \lambda_{i,k}^g v_i^k +\sum_{i\in J} \lambda_{i,k}^h \nabla h_i(x^k)
+\sum_{i=1}^l \nabla \Phi_i(x^k)^T \eta_i^k,
\end{eqnarray} 
but the set of vectors
\begin{eqnarray*}
 \{v_i^{k,s}\}_{i\in \mathcal{I}^k}
 \cup \{\nabla h_i(y^{k,s})\}_{ i\in J}\cup \bigcup_{\beta_i^k\in A_{i,k}^I\cup A_{i,k}^E, i=1,\cdots,l}\{\nabla \Phi_i(y^{k,s})^T \beta_i^k\}
\end{eqnarray*}
is linearly independent  for some sequences  $y^{k,s}\to x^k$, $y^{k,s}\neq x^k$, $\partial g_i(y^{k,s}) \ni v_i^{k,s}\to v_i^k \  (i\in \mathcal{I}^k)$ as $s\to\infty$,
 any set of 
  linearly independent 
 vectors $A_{i,k}:=A_{i,k}^I\cup A_{i,k}^E$ where $A_{i,k}^I \subseteq A_{\Gamma_i}^I(\Phi_i(x^k))$, $A_{i,k}^E \subseteq A_{\Gamma_i}^E(\Phi_i(x^k))$ satisfying  $0\not =\eta_i^k\in {\cal G}{(A_{i,k}^I,A_{i,k}^E)} \subseteq \mathcal{N}_{{\Gamma_i}}(\Phi_i(x^k))$ and $A_{i,k}:=A_{i,k}^I\cup A_{i,k}^E=\emptyset$ if $\eta_i^k=0$.

By  Remark \ref{nor-lim}, for sufficiently large $k$, we have that ${\cal G}{(A_{i,k}^I,A_{i,k}^E)} \subseteq \mathcal{N}_{{\Gamma_i}}(\Phi_i(x^k))\subseteq \mathcal{N}_{{\Gamma_i}}(\Phi_i(\bar x))$,
$A_{i,k}^I \subseteq A_{\Gamma_i}^I(\Phi_i(\bar x))$, and  $A_{i,k}^E \subseteq A_{\Gamma_i}^E(\Phi_i(\bar x))$,
 $i=1,\cdots, l$.
Since $I_g(x^k)\subseteq \bar{I}_g$ and the set $\bar{I}_g$ and the generator set are both finite, taking a subsequence if necessary, we  may assume that
 $\mathcal{I}^k\equiv I\subseteq \bar{I}_g$ and 
 $A_{i,k}^I\equiv A_i^I$, $A_{i,k}^E\equiv A_i^E$,  $i=1,\cdots, l$.

By the diagonalization law,  there exists a sequence $\{z^k\}$ converging to $\bar x$ such that for each $k$, $\bar{v}_i^{k}\in \partial g_i(z^k)$, $\bar{v}_i^{k}\to{v}_i^{k}, i=1,\cdots,n$ and
\begin{eqnarray}\label{lin0}
 \{\bar{v}_i^{k}\}_{i\in I}
 \cup \{\nabla h_i(z^k)\}_{ i\in J}\cup \bigcup_{\beta_i^k\in A_i^I\cup A_i^E, i=1,\cdots,l}\{\nabla \Phi_i(z^k)^T \beta_i^k\}
\end{eqnarray}
is linearly independent, where $A_i^I\cup A_i^E$ is linearly independent, 
$A_{i}^I \subseteq A_{\Gamma_i}^I(\Phi_i(\bar x))$, $A_{i}^E \subseteq A_{\Gamma_i}^E(\Phi_i(\bar x))$, ${\cal G}{(A_{i}^I,A_{i}^E)} \subseteq  \mathcal{N}_{{\Gamma_i}}(\Phi_i(\bar x))$.

Without loss of generality, assume there exist $(\lambda^g,\lambda^h)$, $\bar{\eta}_i\in {\cal G} (A_{i}^I,A_{i}^E)\subseteq\mathcal{N}_{{\Gamma_i}}(\Phi_i(\bar x))$ and $v_i\in \partial g_i(\bar x)$, $i\in I$ satisfying $\lim_{k\to \infty, k\in K}\frac{(\lambda_k^g,\lambda_k^h)}{\|(\lambda_k^g,\lambda_k^h)\|_2}=(\lambda^g,\lambda^h)$, $\displaystyle\lim_{k\rightarrow \infty,k\in K} \eta_i^k/\|(\lambda^g_k, \lambda_{k}^h)\|_2=\bar{\eta}_i$,  for each $i=1,\cdots, l$ and $v_i^k\to v_i$, $i\in I$.

Dividing by $\|(\lambda^g_k, \lambda_{k}^h)\|_2$ on both sides of $(\ref{persisnew})$ and letting $k\to \infty, k\in K$, 
we have that $(\ref{sum1})$ holds.
From  RCPLD (ii),  
{the set of vectors 
\begin{eqnarray*}
 \{\bar{v}_i^{k}\}_{i\in I}
 \cup \{\nabla h_i(z^k)\}_{ i\in J}\cup \bigcup_{\beta_i^k\in A_i^I\cup A_i^E, i\in L}\{\nabla \Phi_i(z^k)^T \beta_i^k\}
\end{eqnarray*}
is linearly dependent, where $L:=\{i=1,\cdots,l: \bar{\eta}_i\neq 0\}$.}
This implies that the set of vectors in (\ref{lin0}) is linearly dependent. This is a contradiction and thus the proof is completed.
\BOX

 Xu and Ye \cite[Definition 1.1]{xy2} introduced a concept of  RCPLD for the  system
\begin{equation}
g(x)\leq 0, h(x)=0, x\in C, \label{NLP}
\end{equation} where $C $ is a closed set.
Let us call  the  system (\ref{feasibility})  a general system provided that one of  $\Gamma_i$ is a closed set instead of a disjunctive set.   By introducing an auxiliary variable $y$ and considering a similar transformation as in the proof of Theorem \ref{OC2} and applying  Xu and Ye \cite[Definition 1.1]{xy2} to the transformed system  one can easily obtain   
a corresponding RCPLD  for the general system.  The RCPLD for the general system turns out to be the same as the one in Definition \ref{nrcpld}  up to and include  (\ref{nna1}) with the condition after  (\ref{nna1}) changed to the following condition.

\noindent $\bullet$ Then for sufficiently large $k$, 
\begin{eqnarray}\label{seqenceLD}
\{v_i^k\}_{i\in I}\cup \{\nabla h_i(x^k)\}_{i\in J}
\cup  \{ \nabla \Phi_i(x^k)^T\eta_i^k\}_{i\in L},
\end{eqnarray}  where $L:=\{i=1,\cdots,l:\eta_i^k\neq 0\}$,
is linearly dependent, 
{for any $\{x^k\}, \{y^k\}, \{v^k\}, \{\eta^k\}$ such that  {$(x^k,y^k) \to (\bar x,  \Phi(\bar x))$,  $x^k \neq \bar x$, 
$v_i^k \in \partial g_i(x^k)$, $v_i^k\rightarrow v_i$, $\eta_i^k\in \mathcal{N}_{\Gamma_i}({y_i^k})$},} 
$\eta^k\to \bar \eta$ as $k\to\infty$.

In particular taking $y_i^k=\Phi_i(\bar x)$ and $\eta^k=\bar \eta$,  
the above condition demands that 
the family of vectors
\begin{equation}\label{general}
\{v_i^k\}_{i\in I}\cup \{\nabla h_i(x^k)\}_{i\in J}
\cup \bigcup_{i\in L}  \{ \nabla \Phi_i(x^k)^T \bar \eta_i \}
\end{equation}
{where $L:=\{i=1,\cdots,l:\bar{\eta}_i\neq 0\}$,} must be linearly dependent.  Now suppose all $\Gamma_i$ are disjunctive,   
$A^i=A_i^I\cup A_i^E$ is linearly independent, 
 $A_i^I \subseteq A_{\Gamma_i}^I(\Phi_i(\bar x))$, $A_i^E \subseteq A_{\Gamma_i}^E(\Phi_i(\bar x))$, 
$0\not =\bar \eta_i\in {\cal G}{(A_i^I,A_i^E)} \subseteq \mathcal{N}_{{\Gamma_i}}(\Phi_i(\bar x))$ 
{and $A^i=A_i^I\cup A_i^E=\emptyset$ if $\bar \eta_i=0$}.  Then the linear dependence of (\ref{general}) would imply the linear dependence of 
\begin{eqnarray*}
\{v_i^k\}_{i\in I}\cup\{\nabla h_i(x^k)\}_{i\in J}\cup  
\bigcup_{\beta_i\in A_i^I\cup A_i^E,i\in \{1,\cdots,l \}} \{\nabla \Phi_i(x^k)^T \beta_i\}. 
\end{eqnarray*} 
Hence in general MPDC-RCPLD is weaker than 
the RCPLD \cite[Definition 1.1]{xy2} for the general system applied to the disjunctive system.

{Therefore, sufficient conditions for  RCPLD defined in \cite[Definition 1.1]{xy2} such as the Linear Constraint Qualification (i.e., when $g,h,\Phi_i$ are affine mappings), NNAMCQ,  relaxed constant rank constraint qualification (RCRCQ)   introduced in the Section 4 of \cite{xy2} applied to the disjunctive system also imply the MPDC-RCPLD.}

\section{Piecewise RCPLD as  a sufficient condition for the
error bound condition}\label{section4}

In this section,  we will propose an RCPLD-type condition 
 that is sufficient for the error bound property.
Although all our results can be stated for the more general case (\ref{feasibility}) where the constrained set $\Gamma_i$ can be different for each $i$, for simplicity of notation, we consider the following disjunctive system:
\begin{equation}
 g(x)\leq 0,
 h(x)=0,
\Phi_i(x)\in {\Gamma},  i=1,\cdots,l,
 \label{s-system} 
\end{equation}
where $g,h,\Phi$ are defined as in  system (\ref{feasibility}),  $\Gamma:=\bigcup_{r=1}^R C_r\subseteq \mathbb{R}^p$ with $C_r$ being a convex polyhedral set, $r=1,\cdots,R$.
Let sets $P_1,\cdots, P_R$ be a partition  of $\{1,\cdots, l\}$ and denote it by $P:=\{ P_1,\cdots, P_R\}$.  Consider
 the subsystem for the partition $P$:
\begin{eqnarray}
{ \left \{\begin{array}{ll }
 g(x)\leq 0,\ h(x)=0, & \\
\Phi_i (x) \in C_1,\  i\in P_1,& \\
~~~~~~ \vdots &  \\
\Phi_{i}(x) \in C_R,\  i\in P_R. &
  \end{array} \right. 
 \label{subsystem}}
\end{eqnarray}
We denote the set of solutions to the subsystem (\ref{subsystem}) by $\widetilde{\cal F}_P$ and the solution set of (\ref{s-system}) by $\widetilde{\cal F}$. 
Since sets $\{ P_1,\cdots, P_R\}$ form  one of the possible partitions of  $\{1,\cdots, l\}$,  the feasible region $\widetilde{\cal F}_P\subseteq \widetilde{\cal F}$ 
{ and for any $\bar x\in {\widetilde{\cal F}}$, there at least exists one partition $P$ such that $\bar x\in \widetilde{\cal F}_P$.}

\begin{thm}\label{perr}
Let $\bar x\in {\widetilde{\cal F}}$.  
{For any partition of $\{1,\cdots, l\}$ into sets $P_1,\cdots, P_R$ such that $\bar x\in \widetilde{\cal F}_P$,  suppose the error bound property for subsystem  (\ref{subsystem}) holds  at $\bar x$,} i.e., there exist $\kappa_P>0$ and $\varepsilon_P>0$ such that 
\begin{eqnarray*}\label{rerr-1}
d_{\widetilde{\cal F}_P}(x) \leq \kappa_P\left  (\|g^+(x)\| + \|h(x)\|+ \sum_{r\in L(x)} \sum_{i\in P_r}  d_{C_r}(\Phi_{i}(x))\right ),
 \quad \forall x\in \mathbb{B}_{\varepsilon_P}(\bar x),
\end{eqnarray*}
{where $L(x):=\{r=1,\cdots,R:P_r\neq \emptyset\}$.}
 Then $\bar x$ satisfies the error bound property  for (\ref{s-system}).
\end{thm} 
{\bf Proof.} 
Let $\varepsilon\leq\min \varepsilon_P$ and $\kappa:=\max \kappa_P$. Then $\mathbb{B}_{\varepsilon}(\bar x)\subseteq \mathbb{B}_{\varepsilon_P}(\bar x)$.
For any $x\in \mathbb{B}_{\varepsilon}(\bar x)$, from the continuity of $\Phi_i(\cdot)$, there must exist a partition $P^x:=\{ P_1^x,\cdots, P_R^x\}$ of $\{1,\cdots, l\}$, where
\begin{eqnarray*}
&&P_1^x:=\big  \{i\in \{1,\cdots, l\}| d_{\Gamma}(\Phi_i(x))=d_{C_1}(\Phi_i(x)) \big\},\\
&&P_r^x:=\left \{i\in \{1,\cdots, l\}| d_{\Gamma}(\Phi_i(x))=d_{C_r}(\Phi_i(x)), i\notin \cup_{j=1}^{r-1} P_j^x \right \},\ r=2,\cdots, R.
\end{eqnarray*}
Then from the error bound property for subsystem  (\ref{subsystem}), for any $x\in \mathbb{B}_{\varepsilon}(\bar x)$, we have that
\begin{eqnarray*}
d_{\widetilde{\cal F}}(x)&=& \inf \left\{\|x-z\|: z\in {\widetilde{\cal F}} \right \}\leq \inf \left \{\|x-z\|: z\in \widetilde{\cal F}_{P^x} \right \}
=d_{\widetilde{\cal F}_{P^x}}(x)\\
&\leq& \kappa(\|g^+(x)\| + \|h(x)\|+ \sum_{r\in L(x)}\sum_{i\in P_r^x} d_{C_r}(\Phi_i(x)))\\
&=&\kappa (\|g^+(x)\| + \|h(x)\| +  \sum_{i=1}^l d_{\Gamma}(\Phi_i(x))).
\end{eqnarray*}
This completes the proof.
\BOX

Theorem \ref{perr} shows that if every subsystem of the disjunctive system satisfies the local error bound, then the local error bound holds for the whole disjunctive system. Based on this observation, we propose the following piecewise RCPLD for the disjunctive system, where {MPDC-RCPLD} is assumed to be satisfied for any subsystem, to ensure the error bound property.

\begin{defn}[{MPDC-PRCPLD}] \label{pnrcpld}   
 The piecewise  RCPLD ({MPDC-PRCPLD}) is said to hold for system  (\ref{s-system}) at $\bar x\in {\widetilde{\cal F}}$, if MPDC-RCPLD holds for subsystem  (\ref{subsystem})
{for any partition of $\{1,\cdots, l\}$ into sets $P_1,\cdots, P_R$ such that $\bar x\in \widetilde{\cal F}_P$.}
\end{defn}


In Xu and Ye \cite[Theorem 3.2]{xy2}, it was shown that if the RCPLD introduced in   \cite[Definition 1.1]{xy2} holds and $g$ is subdifferentially regular, $h$ is continuously differentiable and $\Gamma$ is Clarke regular, then the local error bound holds for  system (\ref{NLP}). Similarly we can show the following error bound result for disjunctive system (\ref{s-system}).

Before proving the error bound property, we need the following two conclusions first.
\begin{lemma}\label{axillary}  
{Assume that  {MPDC-RCPLD} holds 
for the system 
$$\mathcal{F}_0: g(x)\leq 0, h(x)=0,\Phi(x) \in C$$
 at $\bar x\in \mathcal{F}_0$, $C\subseteq \mathbb{R}^p$ is a convex polyhedral set, then {MPDC-RCPLD} holds at $(\bar x,\bar y)$ with $\bar y=\Phi(\bar x)$ for the system $\chi$:}
\begin{eqnarray*}
  g(x)\leq 0,
 h(x)=0, \Psi(x,y):=\Phi(x)-y=0,\  y\in C.
\end{eqnarray*} 
\end{lemma} 
{\bf Proof.}  (i) It is easy to see that for any $(x,y)$,  
$
\nabla\Psi(x,y)=[\nabla \Phi(x),-I_{p}]\in \mathbb{R}^{p\times (d+p)}
$, where $I_p$ denotes the identity matrix of size $p$.  
It is well-known that the rank of a matrix is equal to its row rank, i.e, the rank of the set of its rows. 
Since $\nabla\Psi(x,y)$ has full row ranks for any $(x,y)$ and $\{\nabla h_i(x)\}_{i=1}^m$ have the same rank for all $x$ sufficiently close to $\bar x$, then 
$$
\left [\begin{array}{cc}
\nabla h(x)^T & \nabla \Phi(x)^T \\
0_{q\times m} & -I_{q} 
\end{array}\right ]
$$
have the same rank  for all $(x,y)$ sufficiently close to $(\bar x, \bar y)$. 
Then RCPLD (i) in Definition \ref{nrcpld} holds for system $\chi$.

(ii) Let $J \subseteq \{1,\cdots,m\}$ be such that $\{\nabla h_i(\bar x)\}_{i\in J}$ is a basis for  span$\{\nabla h_i(\bar x)\}_{i=1}^m$.  
  Let  $I\subseteq \bar{I}_g$ be arbitrary index set. 
Suppose there exist vectors $\lambda^h, \lambda^g,\lambda^{\Phi},\bar \eta$ not all zero with 
$\lambda_i^g\geq 0$ for $i\in I$ such that  $v_i\in \partial g_i(\bar x)$,   $\bar \eta\in \mathcal{N}_{C}(\bar y)$  and 
\begin{eqnarray}
&& 0 =   \sum_{i\in I} \lambda_i^g v_i\times\{0\}^p 
+\sum_{i\in J} \lambda_i^h \nabla h_i(\bar x)\times\{0\}^p 
+\nabla \Psi(\bar x,\bar y)^T\lambda^{\Phi}
 +\{0\}^d\times \bar \eta.\label{cq1} 
\end{eqnarray} 

Since (\ref{cq1}) implies
$$0 = \sum_{i\in I} \lambda_i^g v_i  +\sum_{i\in J} \lambda_i^h \nabla h_i(\bar x)+  \nabla \Phi(\bar x)^T \bar  \eta,$$
by RCPLD (ii) in Definition \ref{nrcpld}, the set of vectors  $\displaystyle\{v_i^k\}_{i\in I}\cup  
\{\nabla h_i(x^k)\}_{i\in J} \cup 
 \bigcup_{\beta\in A} \{\nabla \Phi(x^k)^T \beta\}$ 
 is linearly dependent for any sequences $x^k\to \bar x$, $x^k\not = \bar x$, $v_i^k\rightarrow v_i$, $v_i^k \in \partial g_i(x^k)$ ($i\in I$) as $k\to \infty$ and any 
 linearly independent 
 set of vectors $A:=A^I\cup A^E$, $A^I \subseteq A_{C}^I(\Phi(\bar x))$, $A^E \subseteq A_{C}^E(\Phi(\bar x))$ satisfying  $0\not =\bar \eta\in {\cal G}{(A^I,A^E)} \subseteq \mathcal{N}_{C}(\Phi(\bar x))$ and $A:=A^I\cup A^E=\emptyset$ if $\bar \eta=0$.

 It follows that  the vectors 
 $$\{ v_i^k \times \{0\}^p) \}_{i\in I}\cup
  \{\nabla h_i(x^k)\times\{0\}^p\}_{i\in J} \cup 
  \{\nabla \Psi_i(x^k,y^k)\}_{i=1}^p\cup \{\{0\}^d\times  \beta\}_{\beta\in A}$$
   is linearly dependent for $(x^k,y^k)\to (\bar x,\Phi(\bar x))$, $(x^k,y^k)\neq (\bar x,\Phi(\bar x))$.
Thus RCPLD (ii) for system $\chi$  holds at  $(\bar x,\bar y)$.
\BOX

\begin{lemma}\label{err-1}  
{Assume that  {MPDC-RCPLD} holds 
for the system 
$$\mathcal{F}_1: g(x)\leq 0, h(x)=0,x \in C$$
 at $\bar x\in \mathcal{F}_1$, 
 $g_i(\cdot), i=1,\cdots,n$ are subdifferentially regular  and $C$ is a polyhedral convex set.
 Then the error bound property  holds at $\bar x$ for $\mathcal{F}_1$.}
\end{lemma} 
{\bf Proof.}  Assume for a contradiction that the error bound property for the constraint system $\mathcal{F}_1$ fails, then there exists a sequence $C\ni x^k\to \bar x$ and for sufficiently large $k$, 
\begin{eqnarray}\label{errcon}
d_{\mathcal{F}_1}(x^k)\geq k (\|g^+(x^k)\|+\|h(x^k)\|).
\end{eqnarray}

Obviously $x^k\notin \mathcal{F}_1$. Let $y^k$ be the projector of $x^k$ to $\mathcal{F}_1$.  Then $d_{\mathcal{F}_1}(x^k)=\|x^k-y^k\|\not =0$ and $\displaystyle \lim_{k\to\infty}y^k=\bar x$. 
For each $k$, $y^k$ is an optimal solution of 
\begin{eqnarray*}
({\rm P}'_k)~~~\min && F^k(x):=\|x-x^k\|_2\\
 {\rm s.t.} && g(x)\leq 0, h(x)=0, x\in C.
\end{eqnarray*}

From the Proposition \ref{persist}, the MPDC-RCPLD persists in a neighborhood of $\bar x$,  thus MPDC-RCPLD holds at $y^k$ for system $\mathcal{F}_1$, for $k$ sufficiently large.  From Theorem \ref{OC2}, $y^k$ is a stationary point of $({\rm P}'_k)$. 
By the optimality condition, 
there exist parameters $\lambda_{i,k}^g\geq 0$, $v_i^k\in \partial g_i(y^k)$ for $i\in I(y^k)$ and $\lambda_{i,k}^h, i=1,\cdots,m$, $\eta^k\in \mathcal{N}_{C}(y^k)$ such that 
\begin{eqnarray}\label{kt}
0= \frac{y^k-x^k}{\|y^k-x^k\|_2}+\sum_{i\in I(y^k)} \lambda_{i,k}^g v_i^k +\sum_{i=1}^m \lambda_{i,k}^h \nabla h_i(y^k)+\eta^k.
\end{eqnarray}
There exists a set of vectors $A_{k}:=A_{k}^I\cup A_{k}^E$ where $A_{k}^I \subseteq A_{C}^I(y^k)$, $A_{k}^E \subseteq A_{C}^E(y^k)$ satisfying  $0\not =\eta^k\in {\cal G}{(A_{k}^I,A_{k}^E)} \subseteq \mathcal{N}_{C}(y^k)$ and $A_{k}:=A_{k}^I\cup A_{k}^E=\emptyset$ if $\eta^k=0$.

Assume that $\{\nabla h_i(\bar x)\}_{i\in \mathcal{I}_1}$ with $\mathcal{I}_1\subseteq \{1,\cdots,m\}$  is a basis for span$\{\nabla h_i(\bar x)\}_{i=1}^m$.
From Lemma \ref{lem3-1}, we obtain
$
\mathcal{I}^k\subseteq I(y^k)
$ and $\widetilde{A}_{k}^I\subseteq A_{k}^I, \widetilde{A}_{k}^E\subseteq A_{k}^E$
with 
$\{\bar{\lambda}_{k}^g, \bar{\lambda}_{k}^h\}$, $\bar{\lambda}_{i,k}^g> 0$ for $i\in \mathcal{I}^k$ 
  such that 
\begin{eqnarray}\label{AKKT30}
0\in \frac{y^k-x^k}{\|y^k-x^k\|_2}+ \sum_{i\in \mathcal{I}^k} \bar{\lambda}_{i,k}^g v_i^k+\sum_{i\in \mathcal{I}_1} \bar{\lambda}_{i,k}^h \nabla h_i(y^k)+
{\cal G}(\widetilde{A}_{k}^I,\widetilde{A}_{k}^E)
\end{eqnarray}
and ${\cal G}(\widetilde{A}_{k}^I,\widetilde{A}_{k}^E)\subseteq {\cal G}(A_{k}^I,A_{k}^E) \subseteq \mathcal{N}_{C}(y^k)$
with bounded multipliers $\{(\bar{\lambda}_{k}^g, \bar{\lambda}_{k}^h)\}$ from Theorem \ref{OC2}.

Since $C$ is a polyhedral and hence regular, for any $\bar{\eta}^k\in {\cal G}(\widetilde{A}_{k}^I,\widetilde{A}_{k}^E)\subseteq \mathcal{N}_{C}(y^k)$, 
we have that $\langle \bar{\eta}^k,x^k-y^k \rangle \leq \frac{1}{4}\|x^k-y^k\|_2$. 
%
The error bound property holds for $\mathcal{F}_1$ at $\bar x$ followed from the proof of \cite[Theorem 3.2]{xy2}.
\BOX

\begin{thm}\label{newerrorb}  
Suppose that  $g$ is subdifferentially regular and $h,\Phi_i$ are smooth. Assume that the MPDC-PRCPLD holds at $\bar x\in {\widetilde{\cal F}}$ for  system (\ref{s-system}),  then $\bar x$ satisfies error bound property.
\end{thm}
{\bf Proof.} From Definition \ref{pnrcpld}, MPDC-RCPLD holds at $\bar x\in \widetilde{\cal F}_P$ for subsystem  (\ref{subsystem}), 
where $P:=\{ P_1,\cdots, P_R\}$ is a partition  of $\{1,\cdots, l\}$.
From Lemma \ref{axillary}, it is not difficult to see that the MPDC-RCPLD holds at $(\bar x,\bar y)$ with $\bar y=\Phi(\bar x)$ for the system $\widetilde{\chi}_P$:
\begin{eqnarray*}
  g(x)\leq 0,
 h(x)=0,\Phi(x)-y=0,\  y_i \in {C_r},  i\in P_r,\ r=1,\cdots, R.
\end{eqnarray*}
From Lemma \ref{err-1}, the error bound holds at $(\bar x,\bar y)\in \widetilde{\chi}_P$.

For any $x$ around $\bar x$,
\begin{eqnarray*}
d_{\widetilde{\cal F}_P}(x)=\inf \left \{\|x-z\|: z\in \widetilde{\cal F}_P \right \}\leq \inf \left \{\|(x,\Phi(x))-(z,\Phi(z))\|: z\in \widetilde{\cal F}_P\right \}
=d_{\widetilde{\chi}_P}(x,\Phi(x)).
\end{eqnarray*}
Thus the error bound property holds for subsystem  (\ref{subsystem})  at $\bar x\in \widetilde{\cal F}_P$ and the proof follows from Theorem \ref{perr}.
\BOX

We now discuss the relationship between the MPDC-PRCPLD and MPDC-RCPLD. In the following theorem,  by using the relationship between the generator sets of the limiting normal cone of $\Gamma$ and the ones for $C_r$ in Lemma \ref{Lem2.1},   we show that if the dimension of the space where $\Gamma$ is contained is less or equal to two, then the MPDC-PRCPLD implies MPDC-RCPLD for system (\ref{s-system}).

\begin{thm}\label{2dim} 
Assume that $p\leq 2$. Then the MPDC-PRCPLD implies MPDC-RCPLD for system  (\ref{s-system}).
\end{thm} 
{\bf Proof.} 
If $R=1$, then $\Gamma$ itself is a convex polyhedral set and the MPDC-PRCPLD coincides with MPDC-RCPLD. Assume that $R\geq 2$.
 Suppose that MPDC-PRCPLD holds at $\bar x$, then MPDC-RCPLD satisfies at $\bar x$ for the subsystem (\ref{subsystem}) for each partition $P_1,\cdots, P_R$ of $\{1,\cdots, l\}$ such that $\bar x\in \widetilde{\cal F}_P$. We only need to prove that MPDC-RCPLD (ii)  holds for  (\ref{s-system}).
Suppose $J\subseteq \{1,\cdots,m\}$ is the index set that  $\{\nabla h_i(\bar x)\}_{i\in {J}}$  is a basis of $ \mbox{span} \{\nabla h_i(\bar x)\}_{i=1}^m.$ 
Suppose that there exist ${I}\subseteq \bar I_g$ and a nonzero vector $(\lambda^h,\lambda^g,\bar{\eta})$, $v_i \in \partial g_i(\bar x)$, $\lambda_i^g\geq 0$ for $i\in {I}$,  $\bar{\eta}_i\in \mathcal{N}_{\Gamma}(\Phi_i(\bar x))$, for $i\in\{1,\cdots, l\}$  such that 
\begin{eqnarray}\label{2dim1}
0 = \sum_{i\in I} \lambda_i^g v_i+\sum_{i\in J} \lambda_i^h \nabla h_i(\bar x)  + \sum_{i=1}^l \nabla \Phi_i(\bar x)^T \bar  \eta_i.
\end{eqnarray}
We want to prove the linearly dependence of
\begin{eqnarray}
\{v_i^k\}_{i\in I}\cup\{\nabla h_i(x^k)\}_{i\in J}\cup  \bigcup_{\beta_i\in A_i^I\cup A_i^E, i\in\{1,\cdots,l\}} \{\nabla \Phi_i(x^k)^T \beta_i\}, \label{normalvectorphinew}
\end{eqnarray} 
 for any 
$\{x^k\}, \{v^k\}$ such that $x^k\rightarrow \bar x$, $x^k\neq \bar x$, $v^k_i\in \partial g_i(x^k)$, $v^k_i\rightarrow v_i$ and any sets $ A_i^I\subseteq A_{\Gamma}^I(\Phi_i(\bar x)),  A_i^E\subseteq A_{\Gamma}^E(\Phi_i(\bar x))$ such that $A_i:=A_i^I\cup A_i^E$ is linearly independent and $0\not =\bar \eta_i\in {\cal G}{(A_i^I,A_i^E)} \subseteq \mathcal{N}_{{\Gamma}}(\Phi_i(\bar x))$ and $A_i=A_i^I\cup A_i^E=\emptyset$ if $\bar \eta_i=0$.

Since $p\leq 2$, $A_i$ contains at most two vectors. 
Now we will construct a certain partition $P_1,\cdots, P_R$ of $\{1,\cdots, l\}$  such that $\bar x\in \widetilde{\cal F}_P$, which means that $P_1\cup P_2\cup \cdots\cup P_R=\{1,\cdots, l\}$ and $P_j\cap P_k=\emptyset$ for $j\neq k$, $j, k\in \{1,\cdots, l\}$.
Consider the following cases:
\begin{itemize}
\item[{\rm (i)}]
For the index $i\in \{1,\cdots, l\}$ such that
 $A_i=\emptyset$, i.e., $\bar \eta_i=0$. 
 
We select an index $r\in \{1,\cdots, R\}$ satisfying  $\Phi_i(\bar x)\in C_r$, setting  $i\in P_r$ and $A_i^r:=\emptyset$.

\item[{\rm (ii)}] 
For the index $i\in \{1,\cdots, l\}$ such that $A_i$ contains only one vector. 

In this case, there exists $\alpha_i$ such that $A_i=\{\alpha_i \bar{\eta}_i\}\subseteq A_{\Gamma}(\Phi_i(\bar x))$.
From (\ref{unionnorm}) and Lemma \ref{Lem2.1}, we can select an index $s\in \{1,\cdots, R\}$ such that  
$A_i\subseteq A_{C_s}(\Phi_i(\bar x))$.
Set $i\in P_s$ and $A_i^s:=A_i$.

\item[{\rm (iii)}] 
For the index $i\in \{1,\cdots, l\}$ such that $A_i$ contains two vectors. 

In this case, $A_i$ is a basis of $\mathbb{R}^2$.
From (\ref{unionnorm}), 
we can select an index $t\in \{1,\cdots, R\}$ such that  
{$\bar{\eta}_i\in {\cal G}{(A_i^I,A_i^E)}\subseteq {\mathcal{N}}_{C_t}(\Phi_i(\bar x))$.
The generator $((A_i^t)^I, (A_i^t)^E)$ of ${\mathcal{N}}_{C_t}(\Phi_i(\bar x))$ must contain two linearly independent vectors.}
Set $i\in P_t$ and 
$A_i^t:=(A_i^t)^I\cup (A_i^t)^E$.
\end{itemize}
 This way we  construct a partition $P_1,\cdots, P_R$ of $\{1,\cdots, l\}$ such that (\ref{2dim1}) reduces to
\begin{eqnarray*}\label{sum}
0 = \sum_{i\in I} \lambda_i^g v_i+\sum_{i\in J} \lambda_i^h \nabla h_i(\bar x)  + \sum_{i\in P_1} \nabla \Phi_i(\bar x)^T \bar  \eta_i+ \sum_{i\in P_2} \nabla \Phi_i(\bar x)^T \bar  \eta_i+\cdots+ \sum_{i\in P_R} \nabla \Phi_i(\bar x)^T \bar  \eta_i,
\end{eqnarray*}
with $\bar  \eta_i \in  \mathcal{N}_{C_{r}}(\Phi_i(\bar x))$, $i\in P_r$.
From MPDC-RCPLD  (ii) for the subsystem (\ref{subsystem}) for this partition, the family of vectors 
\begin{eqnarray}\label{normalvectorphi-1}
\{v_i^k\}_{i\in {I}}\cup\{\nabla h_i(x^k)\}_{i\in {J}}\cup   \bigcup_{r=1,\cdots,R}\bigcup_{\beta_i^r \in A_i^r, i\in P_r} \{\nabla \Phi_{i}(x^k)^T \beta_i^r\}
\end{eqnarray} 
must be linearly dependent, where $A_i^r$ is defined as in (i)-(iii) for $i\in P_r$, $r=1,\cdots, R$. {Since the sets $A_i^r$ and $A_i$ can be linearly expressed by each other,}
 for $i\in P_r$,  the linear dependence of (\ref{normalvectorphi-1}) implies the the linear dependence of the family of vectors (\ref{normalvectorphinew}).
 Therefore  the proof is completed.
\BOX

For $p\geq 3$,  MPDC-PRCPLD may not imply MPDC-RCPLD. 
The reason is that the inclusion $A_\Gamma (\Phi_i(\bar x)) \subseteq \bigcup_{r=1}^R A_{C_r}(\Phi_i(\bar x))$ may not hold, which can be seen from 
the following example.
\begin{example}\label{Counterexample}
Consider the following disjunctive system:
\begin{eqnarray}\label{expd}
&& h_1(x,y,z):= x^2+xy+x+y+z=0,\nonumber\\
&& h_2(x,y,z):= x-3y-2z=0,\\
&&\Phi(x,y,z):=(x^2-y+z,x+3y^2-z ,-x+2y+z^2 )\in \Gamma,\nonumber
\end{eqnarray}
where $\Gamma:=C_1\cup C_2$ with 
$$C_1:=  \mathbb{R}_+ \times \mathbb{R}_-\times \mathbb{R}_+= \{ w\in \mathbb{R}^3: \langle -e_1,  w\rangle\leq 0,  \langle e_2, w\rangle\leq 0,  \langle e_3, w\rangle\leq 0\}$$ and $C_2:= \{ w\in \mathbb{R}^3: \langle a_i, w\rangle\leq 0, i=1,2\}$, $a_1:=(0.5,-0.5,0.5), a_2:=(-0.5,1,-1)$, 
 see Fig. 1.

Consider the point $(\bar x,\bar y,\bar{z})=(0,0,0)$. 
From (\ref{regularcone}), 
the regular normal cone at $(0,0,0)$ is 
{$\widehat{\mathcal{N}}_{\Gamma}(0,0,0)={\mathcal{N}}_{C_1}(0,0,0)\cap {\mathcal{N}}_{C_2}(0,0,0)={\rm cone}\{-e_1,e_2, -e_3\} \cap {\rm cone}\{a_1,a_2\}={\rm cone}\{a_2,a_3\}$, where $a_3:=(0,1,-1)$.}
\begin{figure}
\begin{minipage}[t]{0.5\linewidth}
  \centering
    \includegraphics[width=1.0\textwidth]{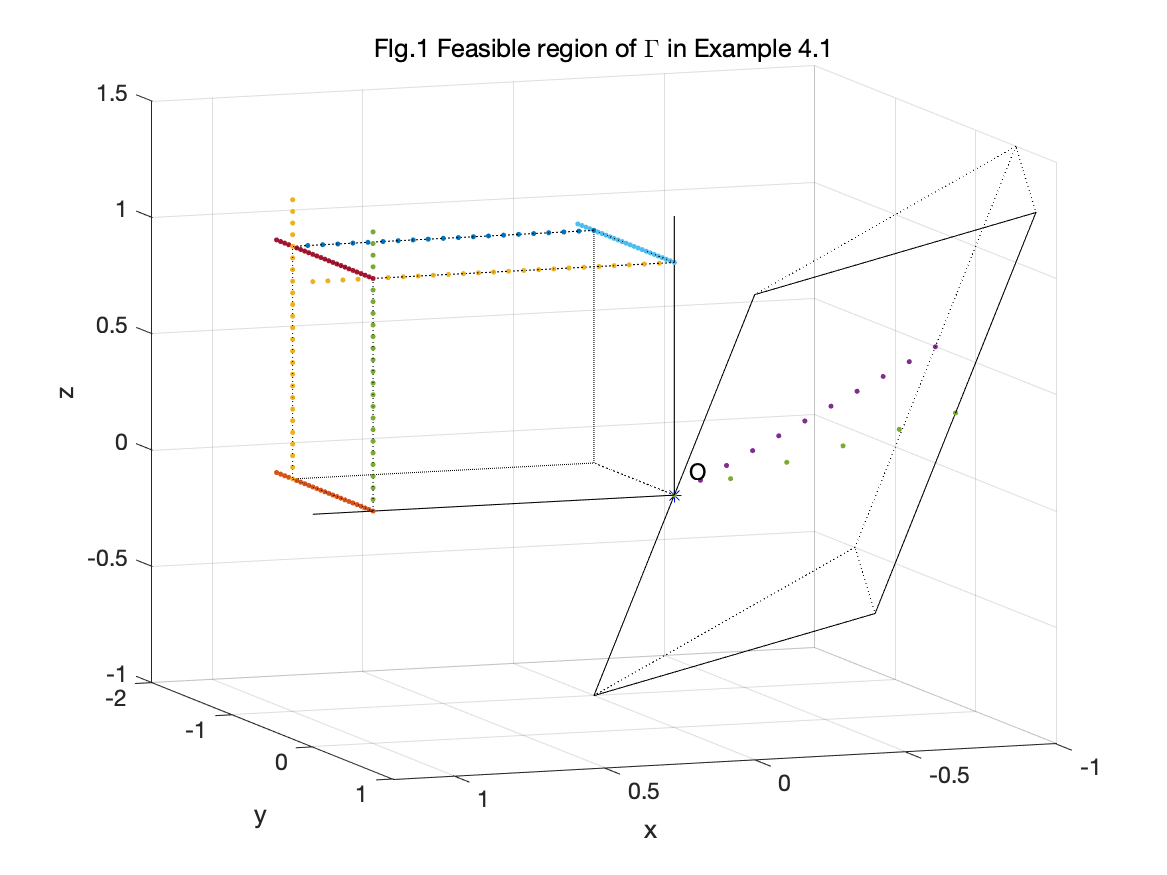}
    \end{minipage}%
\hfill
\begin{minipage}[t]{0.5\linewidth}
\centering
        \includegraphics[width=1.0\textwidth]{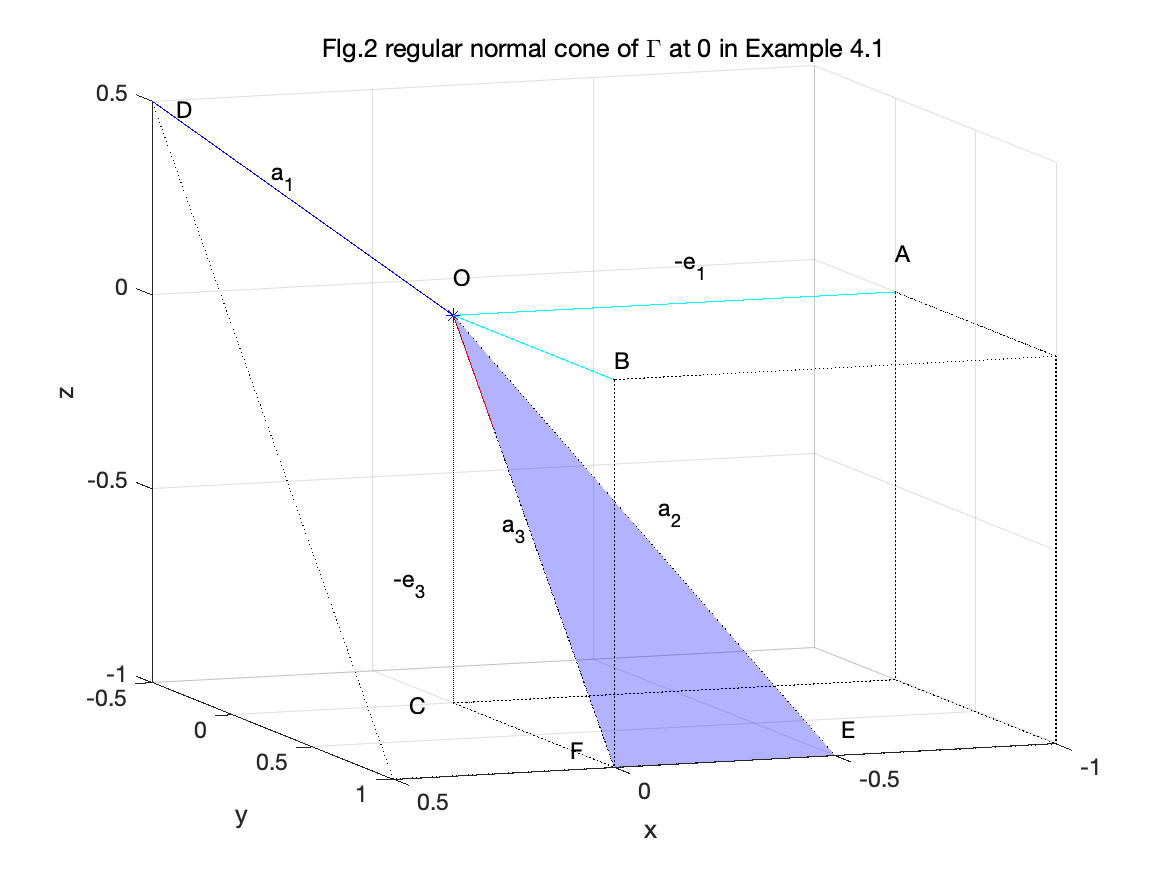}
\end{minipage}
\end{figure}
Similarly we can calculate the regular normal cone of $\Gamma$ at other points and  obtain
\begin{eqnarray*}
\displaystyle \widehat{ \mathcal{N}}_{\Gamma}(x,y,z)=\left \{ \begin{array}{ll}
{\rm cone}\{a_2,a_3\} & \mbox{ if } x=y=z=0,\\
{\rm cone} \{a_1\}  &\mbox{ if } \langle a_1, (x,y,z) \rangle =0, \langle a_2,  (x,y,z)\rangle <0\\
{\rm cone}\{a_2\}& \mbox{ if } \langle a_2,(x,y,z)\rangle  =0,  \langle a_1, (x,y,z)\rangle<0\\
{\rm cone}\{a_1,a_2\}&\mbox{ if }  \langle a_1, (x,y,z)\rangle = \langle a_2, (x,y,z)\rangle =0, (x,y,z)\neq (0,0,0)\\
{\rm cone}\{e_2,-e_3\}&\mbox{ if } x> 0, y=z=0,\\
{\rm cone}\{-e_1,-e_3\}&\mbox{ if } y< 0, x=z=0,\\
{\rm cone}\{-e_1,e_2\}&\mbox{ if } z> 0,  x=y=0,\\
{\rm cone}\{e_2\}&\mbox{ if } y=0, x> 0, z> 0,\\
{\rm cone}\{-e_1\}&\mbox{ if } x=0, y< 0, z> 0,\\
{\rm cone}\{-e_3\}&\mbox{ if } z=0, x> 0, y< 0.\\
\end{array} \right .
\end{eqnarray*}
It is easy to get $\widehat{A}_{\Gamma}(x,y,z)$ from the above representation of regular normal cones.
From (\ref{limitingcone}),   the limiting normal cone at $(0,0,0)$ is
\begin{eqnarray*}
{\mathcal{N}}_{\Gamma}(0,0,0)=\bigcup_{(x,y,z)\in \mathbb{R}^3} \widehat{ \mathcal{N}}_{\Gamma}(x,y,z)= {\rm cone} \{a_1,a_2\}\cup {\rm cone}\{e_2, -e_3\}\cup {\rm cone}\{-e_1,e_2\}\cup {\rm cone} \{-e_1,-e_3\}
\end{eqnarray*}
and from (\ref{Ax}), 
$A_{\Gamma}(0,0,0)=\bigcup_{(x,y,z)\in \mathbb{R}^3} \widehat{A}_{\Gamma}(x,y,z)=\{a_1,a_2,a_3,-e_1,e_2, -e_3\}.$ 
It is obvious that $A_{\Gamma}(0,0,0)\supset A_{C_1}(0,0,0)\cup A_{C_2}(0,0,0)=\{a_1,a_2,-e_1,e_2, -e_3\}$, which {contradicts $A_\Gamma (\Phi(\bar x)) \subseteq \bigcup_{r=1}^2 A_{C_r}(\Phi(\bar x))$.}
%

From {an} easy calculation,
\begin{eqnarray*}
&&\nabla h_1(\bar x,\bar y,\bar z)=
\left(\begin{array}{c}
2\bar x+\bar y+1 \\
\bar x+1\\
1
\end{array} \right)
=\left(\begin{array}{c}
1 \\
1\\
1
\end{array} \right),\ 
\nabla h_2(\bar x,\bar y,\bar z)=
\left(\begin{array}{c}
1 \\
-3\\
-2 
\end{array} \right),\\ 
&&\nabla \Phi(\bar x,\bar y,\bar z)=
\left(\begin{array}{ccc}
2\bar x&-1&1 \\
1&6\bar y& -1 \\
-1&2&2\bar z
\end{array} \right)
=\left(\begin{array}{ccc}
0&-1&1 \\
1&0& -1 \\
-1&2&0
\end{array} \right).
\end{eqnarray*}

We first show that  MPDC-RCPLD fails at $(0,0,0)$. Set $A^I=\{a_3\}\subseteq A_{\Gamma}^I(0,0,0)$, $\bar{\eta}=a_3\in \mathcal{G}(A^I,\emptyset)\subseteq {\mathcal{N}}_{\Gamma}(0,0,0)$, from easy calculation, we have that 
$ \nabla \Phi(\bar x,\bar y,\bar{z})^T \bar \eta=(2,-2,-1)$ and
\begin{eqnarray*}
0 =- \nabla h_1(\bar x,\bar y,\bar{z}) - \nabla h_2(\bar x,\bar y,\bar{z})+ \nabla \Phi(\bar x,\bar y,\bar{z})^T \bar \eta,
\end{eqnarray*}
but for a  sequence $(x^k,y^k,z^k):=(0,0,-\frac{1}{2k})$,  the set of vectors
$$\{\nabla h_1(x^k,y^k,z^k), \nabla h_2(x^k,y^k,z^k), \nabla \Phi(x^k,y^k,z^k)^T a_3\}$$
is linearly independent for any $k\neq 0$.  Thus  MPDC-RCPLD fails at $(0,0,0)$.

We now consider the two subsystems 
\begin{eqnarray*}
\mbox{subsystem 1 } &&h_1(x,y,z)=0, h_2(x,y,z)=0, \Phi(x,y,z)\in C_1.\\
\mbox{subsystem 2 }  &&h_1(x,y,z)=0, h_2(x,y,z)=0, \Phi(x,y,z)\in C_2.
\end{eqnarray*}
For subsystem 1, assume there exists a nonzero vector $(\mu_1, \mu_2,\eta)$, $\eta:=(\eta_1,\eta_2,\eta_3)\in {\mathcal{N}}_{C_1}(0,0,0)$, i.e., $\eta_i\leq 0$, $i=1,3$, $\eta_2\geq 0$
such that 
\begin{eqnarray*}
0 &=& \mu_1 \nabla h_1(\bar x,\bar y,\bar{z})+ \mu_2 \nabla h_2(\bar x,\bar y,\bar{z})+  \nabla \Phi(\bar x,\bar y,\bar{z})^T \eta
\\& =&
\mu_1 \left(\begin{array}{c}
1 \\
1\\
1
\end{array} \right)
+ \mu_2 \left(\begin{array}{c}
1 \\
-3\\
-2
\end{array} \right)
+ \eta_1\left(\begin{array}{c}
0 \\
-1\\
1
\end{array} \right)
+ \eta_2\left(\begin{array}{c}
1 \\
0\\
-1
\end{array} \right)
+ \eta_3\left(\begin{array}{c}
-1 \\
2\\
0
\end{array} \right)
\end{eqnarray*}
If two of $\eta_1,\eta_2,\eta_3$ are equal to zero, without loss of generality we assume $\eta_2=\eta_3=0$. Since the vectors $(1,1,1), (1,-3,-2), (0,-1,1)$ are linearly independent, we must have that $\mu_1=\mu_2=0, \eta_1=0$, which is a contradiction. Similarly, we can prove that
for any $\mu_1, \mu_2\in \mathbb{R}$, the above equation holds provided that at least two of $\eta_1,\eta_2,\eta_3$ do not equal to zero. We assume that $\eta_1\neq 0,\eta_2\neq 0$, then for any sequence $(x^k,y^k,z^k)$, 
$$\{\nabla h_1(x^k,y^k,z^k), \nabla h_2(x^k,y^k,z^k), \nabla \Phi(x^k,y^k,z^k)^T (-e_1), \nabla \Phi(x^k,y^k,z^k)^T e_2\}$$
 must linearly dependent since any four vectors in space $\mathbb{R}^3$ are linearly dependent. Thus MPDC-RCPLD holds for subsystem 1.

For subsystem 2, assume there exists a nonzero vector $(\mu'_1, \mu'_2,\eta')$,   $\eta':=\eta'_1 a_1+\eta'_2 a_2\in cone\{a_1,a_2\}= {\mathcal{N}}_{C_1}(0,0,0)$, $\eta'_i\geq 0$ such that
\begin{eqnarray*}
0 &=& \mu'_1 \nabla h_1(\bar x,\bar y,\bar{z})+ \mu'_2 \nabla h_2(\bar x,\bar y,\bar{z})+  \nabla \Phi(\bar x,\bar y,\bar{z})^T \eta'\\
&=&\mu'_1 \left(\begin{array}{c}
1 \\
1\\
1
\end{array} \right)
+ \mu'_2 \left(\begin{array}{c}
1 \\
-3\\
-2
\end{array} \right)
+ \eta'_1\left(\begin{array}{c}
-1 \\
0.5\\
1
\end{array} \right)
+ \eta'_2\left(\begin{array}{c}
2\\
-1.5\\
-1.5
\end{array} \right)
\end{eqnarray*}
Similarly as subsystem 1, $\eta'_1\neq 0,\eta'_2\neq 0$ and
{ the MPDC-RCPLD holds for subsystem 2. Therefore the MPDC-PRCPLD holds at $(0,0,0)$ for the system (\ref{expd}).}
\end{example}

Now we recall and  introduce some MPDC variants of constraint qualifications that are stronger than MPDC-RCPLD and study their relationships with MPDC-RCPLD. Note that there are other constraint qualifications such as the generalized Guignard CQ, generalized Abadie CQ and the (directional) PQ-/ pseudo-/quasi-normality conditions for MPDC; see Benko et al. \cite{bch}, Flegel et al. \cite{fko}. But since the generalized Guignard CQ, generalized Abadie CQ 
are defined in terms of tangent cones and normal cones of the feasible region and hence not verifiable
and the (directional) pseudo-/quasi-normality conditions are not directly related to RCPLD for MPDC, we do not review them here. 

\begin{defn}\label{cqmd}
 Let $\bar{x}\in {\cal F}$.
The {\em MPDC-LICQ} holds at $\bar x$ if the following condition holds (Mehlitz \cite[Definition 3.1]{mmpscnew}): Let $\Gamma_i:=\bigcup_{r=1}^{R_i}C_r^i$, $C_r^i$ is a convex polyhedral  set,   
$I_{\Gamma_i}(\bar x):=\{r\in\{1,\cdots,R_i\}: \Phi_i(\bar x)\in C_r^i\}$, $i=1,\cdots,l$.  There exists no nonzero vector $(\lambda^h,\lambda^g,\bar \eta) \in \mathbb{R}^{m}\times \mathbb{R}^{n}\times \Pi_{i=1}^l \mathbb{R}^{ p_i}$ with 
$v_i \in \partial g_i(\bar x), i\in \bar{I}_g$ and $\displaystyle \bar \eta_i\in \sum_{r\in I_{\Gamma_i}(\bar x)} {\rm span }\mathcal{N}_{C_r^i}(\Phi_i(\bar x))$  such that
\begin{eqnarray}
 0 = \sum_{i\in \bar{I}_g} \lambda_i^g v_i+\sum_{i=1}^m \lambda_i^h \nabla h_i(\bar x)  + \sum_{i=1}^l \nabla \Phi_i(\bar x)^T \bar  \eta_i.\label{nnal}
\end{eqnarray}  

The {\em MPDC-NNAMCQ} holds at $\bar x$ if for any $v_i \in \partial g_i(\bar x), i\in \bar{I}_g, \bar \eta_i\in \mathcal{N}_{{\Gamma_i}}(\Phi_i(\bar x)), i=1,\cdots,l $, there exists no nonzero vector $(\lambda^h,\lambda^g,\bar \eta) \in \mathbb{R}^{m}\times \mathbb{R}^{n}\times \Pi_{i=1}^l \mathbb{R}^{ p_i}$ with $\lambda_i^g\geq 0$ such that $(\ref{nnal})$ holds.


The {\em MPDC-CRCQ} holds at $\bar x$  if for every index sets $J \subseteq \{1,\cdots,m\}$, $I\subseteq \bar{I}_g$,  $L\subseteq \{1,\cdots,l\}$, every $v_i \in \partial g_i(\bar x), i\in I$, $\bar \eta_i\in \mathcal{N}_{{\Gamma_i}}(\Phi_i(\bar x)), i\in L $, the ranks of
$$\{v_i\}_{i\in I}\cup \{\nabla h_i(\bar{x})\}_{i\in J}\cup \bigcup_{\beta_i\in A_i^I\cup A_i^E, i\in L} \{\nabla \Phi_i(\bar{x})^T \beta_i\} 
$$ and 
$$\{v_i^k\}_{i\in I}\cup \{\nabla h_i(x^k)\}_{i\in J}\cup \bigcup_{\beta_i\in A_i^I\cup A_i^E, i\in L} \{\nabla \Phi_i(x^k)^T \beta_i\} 
$$
are the same, for any 
$\{x^k\}, \{v^k\}$ such that  $x^k\rightarrow \bar x$, $x^k\neq \bar x$, 
 $v^k_i\in \partial g_i(x^k)$, $v^k_i\rightarrow v_i$,
any set of  linearly independent vectors $A_i:=A_i^I\cup A_i^E$ where $A_i^I \subseteq A_{\Gamma_i}^I(\Phi_i(\bar x))$, $A_i^E \subseteq A_{\Gamma_i}^E(\Phi_i(\bar x))$ satisfying  $0\not =\bar \eta_i\in {\cal G}{(A_i^I,A_i^E)} \subseteq \mathcal{N}_{{\Gamma_i}}(\Phi_i(\bar x))$ and $A_i:=A_i^I\cup A_i^E=\emptyset$ if $\bar \eta_i=0$.

The {\em MPDC-RCRCQ} holds at $\bar x$ if the index set $J$ is taken as $\{1,\cdots, m\}$ in MPDC-CRCQ.

The {\em MPDC-ERCPLD} holds at $\bar x$ if
conditions (i)-(ii) are satisfied.  
\begin{itemize}
\item[{\rm (i)}]
 $\{\nabla h_i(x)\}_{i=1}^m$
has constant rank,  for any $x\in \mathbb{B}(\bar x)$,  a neighborhood of $\bar x$. 
\item[{\rm (ii)}]   Suppose $J\subseteq \{1,\cdots,m\}$ is the index set that 
$\{\nabla h_i(\bar x)\}_{i\in J}$  is a basis of $ \mbox{span} \{\nabla h_i(\bar x)\}_{i=1}^m.$ 
If  there exist $I\subseteq \bar I_g$ and $L\subseteq \{1,\cdots,l\}$, 
 $v_i \in \partial g_i(\bar x), i\in I$, $\bar \eta_i\in \mathcal{N}_{{\Gamma_i}}(\Phi_i(\bar x)), i\in L$ such that 
\begin{eqnarray}\label{poslinen}
( \{v_i, i\in I\},\ \{\nabla h_i(\bar{x}), i\in J;\ \nabla \Phi_i(\bar{x})^T \beta_i, \beta_i\in A_i^I\cup A_i^E, i\in L\})
\end{eqnarray}  
is positive linearly dependent, then  for  sufficiently large $k$,
\begin{eqnarray*}
\{v_i^k\}_{i\in I}\cup\{\nabla h_i(x^k)\}_{i\in J}\cup  \bigcup_{\beta_i\in A_i^I\cup A_i^E, i\in L} \{\nabla \Phi_i(x^k)^T \beta_i\} 
\end{eqnarray*} 
is linearly dependent, for any $\{x^k\}, \{v^k\}$ such that  $x^k\rightarrow \bar x$, $x^k\neq \bar x$, 
 $v^k_i\in \partial g_i(x^k)$, $v^k_i\rightarrow v_i$,
any set of  linearly independent vectors $A_i:=A_i^I\cup A_i^E$ where $A_i^I \subseteq A_{\Gamma_i}^I(\Phi_i(\bar x))$, $A_i^E \subseteq A_{\Gamma_i}^E(\Phi_i(\bar x))$ satisfying  $0\not =\bar \eta_i\in {\cal G}{(A_i^I,A_i^E)} \subseteq \mathcal{N}_{{\Gamma_i}}(\Phi_i(\bar x))$ and $A_i:=A_i^I\cup A_i^E=\emptyset$ if $\bar \eta_i=0$.
\end{itemize}

The {\em MPDC-CPLD} holds at $\bar x$ if the following conditions hold.  Suppose that 
 there exist index sets $J \subseteq \{1,\cdots,m\}$, $I\subseteq \bar{I}_g$ and $L\subseteq \{1,\cdots,l\}$,
 a nonzero vector $(\lambda^h,\lambda^g,\bar \eta) \in \mathbb{R}^{m}\times \mathbb{R}^{n}\times \Pi_{i=1}^l \mathbb{R}^{ p_i}$ with $\lambda_i^g\geq 0,$  $v_i \in \partial g_i(\bar x), i\in I$  and $\bar \eta_i\in \mathcal{N}_{{\Gamma_i}}(\Phi_i(\bar x)), i\in L $ such that  
\begin{eqnarray*}
 0 = \sum_{i\in I} \lambda_i^g v_i+\sum_{i\in J} \lambda_i^h \nabla h_i(\bar x)  + \sum_{i\in L} \nabla \Phi_i(\bar x)^T \bar  \eta_i.
\end{eqnarray*}  
Then  for  sufficiently large $k$,
\begin{eqnarray*}
\{v_i^k\}_{i\in I}\cup\{\nabla h_i(x^k)\}_{i\in J}\cup  \bigcup_{\beta_i\in A_i^I\cup A_i^E, i\in L} \{\nabla \Phi_i(x^k)^T \beta_i\} 
\end{eqnarray*} 
is linearly dependent, for any 
$\{x^k\}, \{v^k\}$ such that  $x^k\rightarrow \bar x$, $x^k\neq \bar x$, 
 $v^k_i\in \partial g_i(x^k)$, $v^k_i\rightarrow v_i$,
any set of  linearly independent vectors $A_i:=A_i^I\cup A_i^E$ where $A_i^I \subseteq A_{\Gamma_i}^I(\Phi_i(\bar x))$, $A_i^E \subseteq A_{\Gamma_i}^E(\Phi_i(\bar x))$ satisfying  $0\not =\bar \eta_i\in {\cal G}{(A_i^I,A_i^E)} \subseteq \mathcal{N}_{{\Gamma_i}}(\Phi_i(\bar x))$ and $A_i:=A_i^I\cup A_i^E=\emptyset$ if $\bar \eta_i=0$.

\end{defn}

On one hand,  MPDC-CPLD does not imply the MPDC-ERCPLD.
{Indeed, if the vectors in (\ref{poslinen}) are  positively linearly dependent, i.e.,
there exist $\lambda_i^g\geq 0, i\in I, \lambda_i^h, i\in J, \lambda_i^{\Phi}, i\in L$ satisfying
$$
0=\sum_{i\in I} \lambda_i^g v_i+\sum_{i\in J} \lambda_i^h \nabla h_i(\bar x)  + \sum_{\beta_i\in A_i^I\cup A_i^E,i\in L}  \lambda_i^{\Phi} \nabla \Phi_i(\bar x)^T \beta_i.
$$
From the definition of ${\cal G} (A_i^I, A_i^E)$, the vector $\bar \eta_i:=\displaystyle\sum_{\beta_i\in A_i^I\cup A_i^E}  \lambda_i^{\Phi}  \beta_i$
 may not necessarily belong to ${\cal G} (A_i^I, A_i^E)$ and thus may not belong to $\mathcal{N}_{{\Gamma_i}}(\Phi_i(\bar x))$.}
On the other hand, MPDC-ERCPLD requires $J\subseteq \{1,\cdots,m\}$ being the index set such that $\{\nabla h_i(\bar x)\}_{i\in J}$  is a basis of $ \mbox{span} \{\nabla h_i(\bar x)\}_{i=1}^m$ while in MPDC-CPLD, $J\subseteq \{1,\cdots,m\}$ is an arbitrary set, thus  MPDC-ERCPLD does not imply  MPDC-CPLD either.
The  condition  MPDC-ERCPLD is stronger than  or equal to  MPDC-RCPLD
because the condition (\ref{nna1}) implies the positive linearly dependence of the vectors in (\ref{poslinen}).
Moreover it is easy to see that MPDC-RCRCQ implies  MPDC-ERCPLD. We summarize the relationships among various constraint qualifications  in Figure 3.

%

\begin{figure}%
\centering
	\scriptsize
		 \tikzstyle{format}=[rectangle,draw,thin,fill=white]
		
		 \tikzstyle{test}=[diamond,aspect=2,draw,thin]
		
		\tikzstyle{point}=[coordinate,on grid,]
\begin{tikzpicture}
        \node[format](MPDC-PRCPLD){MPDC-PRCPLD};
        \node[format,right of=MPDC-PRCPLD,node distance=27mm](error bound){error bound};
        \node[format, below of=MPDC-PRCPLD,node distance=10mm](MPDC-NNAMCQ){MPDC-NNAMCQ};
        
        \node[format,left of=MPDC-NNAMCQ,node distance=30mm](MPDC-LICQ){MPDC-LICQ};
        \node[format,right of=MPDC-NNAMCQ,node distance=30mm](MPDC-CPLD){MPDC-CPLD};
	\node[format,right of=MPDC-CPLD,node distance=30mm](MPDC-RCPLD){MPDC-RCPLD};
       \node[format,right of=MPDC-RCPLD,node distance=27mm](M-stationarity){M-stationarity};
       \node[format, below of=MPDC-NNAMCQ,node distance=10mm](MPDC-CRCQ){MPDC-CRCQ};
       \node[format,right of=MPDC-CRCQ,node distance=30mm](MPDC-RCRCQ){MPDC-RCRCQ};
       \node[format,right of=MPDC-RCRCQ,node distance=30mm](MPDC-ERCPLD){MPDC-ERCPLD};
\draw[->](MPDC-LICQ)--(MPDC-NNAMCQ);
\draw[->](MPDC-ERCPLD)--(MPDC-RCPLD);
\draw[->](MPDC-PRCPLD)--(error bound);
\draw[->](MPDC-LICQ)--(MPDC-PRCPLD);
\draw[->](MPDC-NNAMCQ)--(MPDC-CPLD);
\draw[->](MPDC-CPLD)--(MPDC-RCPLD);
\draw[->](MPDC-RCPLD)--(M-stationarity);
\draw[->](MPDC-LICQ)--(MPDC-CRCQ);
\draw[->](MPDC-CRCQ)--(MPDC-RCRCQ);
\draw[->](MPDC-CRCQ)--(MPDC-CPLD);
\draw[->](MPDC-RCRCQ)--(MPDC-ERCPLD);
\draw[->](MPDC-NNAMCQ)--(error bound);
\end{tikzpicture}

\centering{Fig.3 Relation among constraint qualifications for MPDC}
\end{figure}

%
%
%

\section{Applications to ortho-disjunctive systems}
In this section, we consider the ortho-disjunctive programs:
\begin{equation}
 g(x)\leq 0,
 h(x)=0,
\Phi_i(x)\in {\Gamma},  i=1,\cdots,p,
 \label{s-systemO} 
\end{equation}
where $\Phi_i(x)=(\Phi_i^1(x),\Phi_i^2(x))$ with $\Phi_i^1(x):=G_i(x),\Phi_i^2(x):=H_i(x)$, $\Gamma=C_1\cup C_2$, 
with $C_r=[a_1^r, b_1^r]\times [a_2^r, b_2^r]$, $-\infty\leq a_j^r\leq  b_j^r \leq \infty $.
The system (\ref{s-systemO}) is an equivalent reformulation of the mathematical program with ortho-disjunctive constraints (MPODC)/the ortho-disjunctive program introduced in Benko et al. \cite{bch}.
Denote by {${\cal{F}}_O$} the solution set of the  ortho-disjunctive system $(\ref{s-systemO})$.

To apply our results, we rewrite 
$C_r:=\{y\in \mathbb{R}^2:  \langle -e_j,y\rangle\leq -a_j^r, \langle e_j,y\rangle\leq b_j^r, j=1,2\}$, for $r=1,2$. 
For $y\in \mathbb{R}^2$, define the active  index sets 
$I_r^{0}(y):=\{j=1,2:  y_j=a_j^r=b_j^r\},
I_r^{-}(y):=\{j=1,2:  y_j=a_j^r<b_j^r\}, 
I_r^{+}(y):=\{j=1,2:  a_j^r<y_j=b_j^r \}.
$
Then the normal cone to each convex polyhedral set $C_r$ can be calculated as 
\begin{eqnarray}\label{norml}
\mathcal{N}_{C_r}(y)= 
{\rm span}\{e_j: j\in I_r^{0}(y)\}+{\rm cone} \left (\{ -e_j: j\in I_r^{-}(y)\} \cup \{ e_j: j\in I_r^{+}(y)\}\right ).
\end{eqnarray}
 From (\ref{regularcone}), 
{$\widehat{ \mathcal{N}}_{\Gamma}(y)= \mathcal{N}_{C_1}(y)$ if $y\in C_1\setminus C_2$, $\widehat{ \mathcal{N}}_{\Gamma}(y)= \mathcal{N}_{C_2}(y)$ if $y\in C_2\setminus C_1$ and if $y\in C_1\cap C_2$,}
 $$\widehat{ \mathcal{N}}_{\Gamma}(y)\subseteq \mathcal{N}_{C_1}(y) \cap  \mathcal{N}_{C_2}(y).$$
 From (\ref{Ahat}),  $(\widehat{A}_{\Gamma}^I(y), \widehat{A}_{\Gamma}^E(y))$ is the generator set for $\widehat{ \mathcal{N}}_{\Gamma}(y)$, i.e.,
$\widehat{ \mathcal{N}}_{\Gamma}(y)=
\mathcal{G}(\widehat{A}_{\Gamma}^I(y), \widehat{A}_{\Gamma}^E(y))$.
It is obvious that $(\widehat{A}_{\Gamma}^I(y), \widehat{A}_{\Gamma}^E(y))$ is linearly independent and
\begin{equation}\label{eqn5.3}
\widehat{A}_{\Gamma}^I(y)\cup \widehat{A}_{\Gamma}^E(y) \subseteq  \{-e_1,e_1,-e_2,e_2\}.
\end{equation}

For any $i=1,\cdots, p$, let $\bar y:=(G_i(\bar x),H_i(\bar x))$ and $\bar{\eta}_i:=(\bar{\eta}_i^1,\bar{\eta}_i^2) \in \mathcal{N}_{\Gamma}(\bar y)$. Then from (\ref{limitingcone}),  there exists $y\in \mathbb{B}_{\delta}(\bar y)$ such that ${\eta}_i \in 
\widehat{ \mathcal{N}}_{\Gamma}(y)$ and 
we can find $A_i^I\subseteq\widehat{A}_{\Gamma}^I(y), A_i^E\subseteq\widehat{A}_{\Gamma}^E(y)$
satisfying 
\begin{equation}
\eta_i\in {\cal G} (A_i^I,A_i^E) \subseteq  {\cal G} (\widehat{A}_{\Gamma}^I(y),\widehat{A}_{\Gamma}^E(y))= \widehat{ \mathcal{N}}_{\Gamma}(y)\subseteq \mathcal{N}_{{\Gamma}}(\bar y). \label{eqn5.4}
\end{equation}
From (\ref{eqn5.3})-(\ref{eqn5.4}), it is easy to see that $A_i:=A_i^I\cup A_i^E$ can be taken as
\begin{eqnarray}\label{Ai}
A_i= \left \{ \begin{array}{ll}
\{e_1,e_2\}\mbox{ or } \{-e_1,e_2\}\mbox{ or } \{e_1,-e_2\}\mbox{ or } \{-e_1,-e_2\}  & \mbox{ if } {\eta}_i^1\not =0, {\eta}_i^2\not =0,\\
\{e_1\}\mbox{ or } \{-e_1\}  & \mbox{ if } {\eta}_i^1\not =0, {\eta}_i^2=0,\\
\{e_2\}\mbox{ or } \{-e_2\}  & \mbox{ if }  {\eta}_i^1=0, {\eta}_i^2\not =0,\\
\emptyset & \mbox{ if }  {\eta}_i^1={\eta}_i^2 =0.
\end{array}
\right. 
\end{eqnarray}
We explain the above discussions by the complementarity system. 
\begin{example}\label{exme}
Consider the complementary system,  $\Omega_{E}=C_1\cup C_2$, where $C_1:=\{ y\in\mathbb{R}^2:\langle -e_1, y\rangle\leq  0, \langle e_2, y\rangle=0\}$ and $C_2:=\{y\in\mathbb{R}^2: \langle e_1, y\rangle=0,\langle -e_2, y\rangle\leq  0\}$.
From the above discussions, the regular normal cone 
\begin{eqnarray*}
\displaystyle \widehat{ \mathcal{N}}_{\Omega_{E}}(y)=\left \{ \begin{array}{ll}
 \mathcal{N}_{C_1}(y)
=  \{ \lambda e_2: \lambda \in \mathbb{R} \}  & \mbox{ if } y_1>0, y_2=0,\\
\mathcal{N}_{C_2}(y)=  \{ \lambda e_1: \lambda \in \mathbb{R} \}  &\mbox{ if } y_1=0, y_2>0,\\
\mathcal{N}_{C_1}(y)\cap  \mathcal{N}_{C_2}(y)= \{ -\lambda_1 e_1-\lambda_2 e_2: \lambda_1\geq 0,\lambda_2\geq 0\}& \mbox{ if } y_1=y_2=0.
\end{array} \right . 
\end{eqnarray*}
It follows that   
\begin{eqnarray*}
 \begin{array}{ll}
\widehat{A}_{\Omega_{E}}^I(y_1,y_2)=\emptyset,\ \widehat{A}_{\Omega_{E}}^E(y_1,y_2)=\{e_2\},&\mbox{ if } y_1>0, y_2=0,\\
\widehat{A}_{\Omega_{E}}^I(y_1,y_2)=\emptyset,\ \widehat{A}_{\Omega_{E}}^E(y_1,y_2)=\{e_1\},&\mbox{ if } y_1=0, y_2>0,\\
\widehat{A}_{\Omega_{E}}^I(y_1,y_2)=\{-e_1,-e_2\},\ \widehat{A}_{\Omega_{E}}^E(y_1,y_2)=\emptyset, & \mbox{ if } y_1=y_2=0.\end{array}
  \label{AE}
\end{eqnarray*}
 From (\ref{limitingcone}),   the limiting normal cone 
\begin{eqnarray}
\mathcal{N}_{\Omega_{E}}(0,0) &=& \bigcup_{y\in \mathbb{R}^2} \widehat{ \mathcal{N}}_{\Omega_{E}}(y)=
\left  \{-\lambda_1 e_1-\lambda_2 e_2:  {\rm either }\ \lambda_1,\lambda_2 {>} 0\ {\rm or}\ \lambda_1\lambda_2=0 \right \}.\label{normalcone-EC}
\label{limit2}
\end{eqnarray}
From (\ref{eqn5.4}), for any $({\eta}^1,{\eta}^2) \in \mathcal{N}_{\Omega_{E}}(0,0)$, since $A^I\subseteq \widehat{A}_{\Omega_{E}}^I(y_1,y_2), A^E\subseteq \widehat{A}_{\Omega_{E}}^E(y_1,y_2)$ must satisfy
$$\eta\in {\cal G} (A^I,A^E) \subseteq  {\cal G} (\widehat{A}_{\Omega_{E}}^I(y_1,y_2),\widehat{A}_{\Omega_{E}}^I(y_1,y_2))
=\widehat{\mathcal{N}}_{\Omega_{E}}(y_1,y_2)
\subseteq  \mathcal{N}_{\Omega_{E}}(0,0),
$$
one of the following cases must occur:\\
(i) ${\eta}^1\neq 0$, ${\eta}^2\neq 0$, then $({\eta}^1,{\eta}^2) \in \widehat{\mathcal{N}}_{\Omega_{E}}(0,0)$ and thus $A^I=\{-e_1,-e_2\}$, $A^E=\emptyset$; \\
(ii) ${\eta}^1= 0$, ${\eta}^2\neq 0$, then $({\eta}^1,{\eta}^2) \in \widehat{\mathcal{N}}_{\Omega_{E}}(y_1,0)$, 
$y_1\geq 0$ and thus $A^I=\emptyset$, $A^E=\{e_2\}$ or $A^I=\{-e_2\}$, $A^E=\emptyset$; \\
(iii) ${\eta}^1\neq 0$, ${\eta}^2= 0$, then $({\eta}^1,{\eta}^2) \in \widehat{\mathcal{N}}_{\Omega_{E}}(0,y_2)$, 
$y_2\geq 0$ and thus $A^I=\emptyset$, $A^E=\{e_1\}$ or $A^I=\{-e_1\}$, $A^E=\emptyset$; \\
(iv) ${\eta}^1= 0$, ${\eta}^2= 0$, then $A^I=A^E=\emptyset$.
\end{example}

Now applying  MPDC-RCPLD to the ortho-disjunctive system $(\ref{s-systemO})$ {and using the condition (\ref{Ai})}, we obtain the following condition. 
\begin{defn}[{MPODC-RCPLD}]\label{rcpldxi}
 We say that  the $\rm RCPLD$ for the ortho-disjunctive system (\ref{s-systemO}) holds  at $\bar x\in {\cal{F}}_O$ if conditions (i)-(ii) are satisfied.
\begin{itemize}
\item[{\rm (i)}]
 $\{\nabla h_i(x)\}_{i=1}^m$ has constant rank,  for any $x\in \mathbb{B}(\bar x)$,  a neighborhood of $\bar x$. 
\item[{\rm (ii)}]  Suppose $\mathcal{I}_1\subseteq \{1,\cdots,m\}$ is the index set that 
$\{\nabla h_i(\bar x)\}_{i\in \mathcal{I}_1}$  is a basis of $ \mbox{span} \{\nabla h_i(\bar x)\}_{i=1}^m.$ 
Suppose there exist $\mathcal{I}_2\subseteq \bar I_g$, a nonzero vector $(\lambda^h,\lambda^g,\bar \eta^1,\bar{\eta}^2) \in \mathbb{R}^{m}\times \mathbb{R}^{n}\times \mathbb{R}^{ 2p}$ and $v_i \in \partial g_i(\bar x)$, $\lambda_i^g\geq 0, i\in \mathcal{I}_2$ such that
\begin{eqnarray}
&&0 = \sum_{i\in \mathcal{I}_2} \lambda_i^g v_i+\sum_{i\in \mathcal{I}_1} \lambda_i^h \nabla h_i(\bar x)  + \sum_{i=1}^p (\nabla G_i(\bar x) \bar  \eta_i^1+ \nabla H_i(\bar x) \bar  \eta_i^2),\nonumber\label{xinna}\\
&&( \bar \eta_i^1, \bar \eta_i^2) \in \mathcal{N}_{C_1\cup C_2}(G_i(\bar x),H_i(\bar x)), \quad i=1,\cdots,p.\label{xinna1}
\end{eqnarray}  
Then for  sufficiently large $k$,
\begin{eqnarray*}
\{v_i^k\}_{i\in \mathcal{I}_2}\cup\{\nabla h_i(x^k)\}_{i\in \mathcal{I}_1}\cup  \{\nabla G_i(x^k)\}_{i\in \mathcal{I}_3}\cup\{\nabla H_i(x^k)\}_{i\in \mathcal{I}_4}
 \label{normalvectorxi}
\end{eqnarray*}  
is linearly dependent, for any  
$\{x^k\}, \{v^k\}$ such that $x^k\rightarrow \bar x$, $x^k\neq \bar x$, 
 $v^k_i\in \partial g_i(x^k)$, $v^k_i\rightarrow v_i$ as $k\to\infty$, where $\mathcal{I}_3:=\{i=1,\cdots,p: \bar \eta_i^1\neq 0\}$ and $\mathcal{I}_4:= \{i=1,\cdots,p: \bar \eta_i^2\neq 0\}$.
\end{itemize}
\end{defn}

We now define the piecewise RCPLD for the ortho-disjunctive system $(\ref{s-systemO})$.  Around the point $\bar x$, 
consider a partition  of $\{1,\cdots, p\}$ into sets $P, Q$ and the system 
\begin{eqnarray}
{ \left \{\begin{array}{ll }
 g(x)\leq 0,\ h(x)=0, & \\
\Phi_i (x) \in C_1,\  i\in P,& \\
\Phi_{i}(x) \in C_2,\  i\in Q. &
  \end{array} \right. 
 \label{subsystemo}}
\end{eqnarray}
{Then the MPODC-piecewise RCPLD holds at $\bar x$ if MPODC-RCPLD holds for each subsystem (\ref{subsystemo}), i.e., the condition (\ref{xinna1}) in Definition \ref{rcpldxi} is replaced by $( \bar \eta_i^1, \bar \eta_i^2) \in \mathcal{N}_{C_1}(G_i(\bar x),H_i(\bar x))$ for $i\in P$ and $( \bar \eta_i^1, \bar \eta_i^2) \in \mathcal{N}_{C_2}(G_i(\bar x),H_i(\bar x))$ for $i\in Q$.}
From (\ref{norml}), we can find $A_i^I,A_i^E\subseteq \{e_j, j\in I_r^{0}(\bar{y});\ -e_j, j\in I_r^{-}(\bar{y});\ e_j, j\in I_r^{+}(\bar{y})\}$, where $\bar y:=(G_i(\bar x),H_i(\bar x))$  such that $A_i^I,A_i^E$ are disjoint and $A_i^I\cup A_i^E$ is linearly independence and
$${( \bar \eta_i^1, \bar \eta_i^2)}\in {\cal G} (A_i^I,A_i^E) \subseteq 
\mathcal{N}_{C_r}(G_i(\bar x),H_i(\bar x)), $$
where $r=1$ if $i\in P$ and $r=2$ if $i\in Q$.
Thus  
{$A_i$ also satisfies the condition $(\ref{Ai})$.}
Then the piecewise RCPLD for the ortho-disjunctive system $(\ref{s-systemO})$ reduces to the following definition.

\begin{defn}[{MPODC-PRCPLD}]\label{pnrcpldxi}
The piecewise RCPLD  for the ortho-disjunctive system holds at $\bar x\in {\cal{F}}_O$ for $(\ref{s-systemO})$ if all conditions in MPODC-RCPLD hold with (\ref{xinna1}) replaced by 
\begin{eqnarray*}
( \bar \eta_i^1, \bar \eta_i^2) \in \mathcal{N}_{C_1}(G_i(\bar x),H_i(\bar x))\cup \mathcal{N}_{C_2}(G_i(\bar x),H_i(\bar x)), \quad i=1,\cdots,p.
\end{eqnarray*}
\end{defn}

Because we have (see e.g., Mehlitz \cite[Lemma 2.2]{mmpscnew}),
$\mathcal{N}_{{C_1\cup C_2}}(y_1,y_2)\subseteq \mathcal{N}_{C_1}(y_1,y_2)\cup \mathcal{N}_{C_2}(y_1,y_2),$
it is obvious that the MPODC-PRCPLD is stronger than MPODC-RCPLD. This result is consistent with Theorem \ref{2dim}.
The following error bound result of ortho-disjunctive system $(\ref{s-systemO})$ follows from Theorem \ref{newerrorb}. 
\begin{thm}\label{xierr} Suppose that $g$ is subdifferentially  regular, $h, G_i, H_i$ are smooth.
Assume that  the MPODC-PRCPLD holds at $\bar x\in {\cal{F}}_O$ of system $(\ref{s-systemO})$, then $\bar x$ satisfies the error bound holds.
\end{thm} 

We now apply the constraint qualifications for MPDC introduced in Definition \ref{cqmd} to the ortho-disjunctive system $(\ref{s-systemO})$ and study their relationships with MPODC-RCPLD.  
\begin{defn}\label{cqod}
Let $\bar x\in {\cal{F}}_O$.
We say $\bar{x}$ satisfies the {\em MPODC-LICQ}, if for any $v_i \in \partial g_i(\bar x), i\in \bar{I}_g$ and 
$(\bar \eta_i^1,\bar \eta_i^2) \in {\rm span }\mathcal{N}_{C_1}(G_i(\bar x),H_i(\bar x))+{\rm span } \mathcal{N}_{C_2}(G_i(\bar x),H_i(\bar x))$, $i=1\cdots , p$,
there exists no nonzero vector $(\lambda^h,\lambda^g,\bar \eta^1,\bar{\eta}^2) \in \mathbb{R}^{m}\times \mathbb{R}^{n}\times \mathbb{R}^{ 2p}$ such that
\begin{eqnarray}
 0 = \sum_{i\in \bar{I}_g} \lambda_i^g v_i+\sum_{i=1}^m \lambda_i^h \nabla h_i(\bar x)  + \sum_{i=1}^p (\nabla G_i(\bar x) \bar  \eta_i^1+ \nabla H_i(\bar x) \bar  \eta_i^2).\label{nnal-1}
\end{eqnarray}   

We say $\bar{x}$ satisfies the {\em MPODC-NNAMCQ}, if for any $v_i \in \partial g_i(\bar x), i\in \bar{I}_g$, $(\bar \eta_i^1,\bar \eta_i^2)\in \mathcal{N}_{{C_1\cup C_2}}(G_i(\bar x),H_i(\bar x))$, $i=1,\cdots,p$,
there exists no nonzero vector $(\lambda^h,\lambda^g,\bar \eta^1,\bar{\eta}^2) \in \mathbb{R}^{m}\times \mathbb{R}^{n}\times \mathbb{R}^{ 2p}$ with $\lambda_i^g\geq 0$ such that (\ref{nnal-1}) holds.


We say $\bar{x}$ satisfies the {\em MPODC-CRCQ}, if for every index sets $\mathcal{I}_1 \subseteq \{1,\cdots,m\}$, $\mathcal{I}_2\subseteq \bar{I}_g$,  $L\subseteq \{1,\cdots,p\}$ and 
for any $( \bar \eta_i^1, \bar \eta_i^2)\in \mathcal{N}_{{C_1\cup C_2}}(G_i(\bar x),H_i(\bar x)), i\in L$,
 $v_i \in \partial g_i(\bar x), i\in \mathcal{I}_2$,
the ranks of 
$$\{v_i\}_{i\in \mathcal{I}_2}\cup \{\nabla h_i(\bar{x})\}_{i\in \mathcal{I}_1}\cup \{\nabla G_i(\bar x) \}_{i\in \mathcal{I}_3}\cup \{\nabla H_i(\bar x)\}_{i\in \mathcal{I}_4}
$$ and 
$$\{v_i^k\}_{i\in \mathcal{I}_2}\cup \{\nabla h_i(x^k)\}_{i\in \mathcal{I}_1}\cup \{\nabla G_i(x^k)\}_{i\in \mathcal{I}_3}\cup \{\nabla H_i(x^k)\}_{i\in \mathcal{I}_4}
$$
are the same, for any 
$\{x^k\}, \{v^k\}$ such that $x^k\rightarrow \bar x$, $x^k\neq \bar x$, 
 $v^k_i\in \partial g_i(x^k)$, $v^k_i\rightarrow v_i$,
 where $\mathcal{I}_3:=\{i\in L: \bar \eta_i^1\neq 0\}$, $\mathcal{I}_4:= \{i\in L: \bar \eta_i^2\neq 0\}$. 

We say $\bar{x}$ satisfies the {\em MPODC-RCRCQ}, if $\mathcal{I}_1$ is taken as $\{1,\cdots, m\}$ in MPODC-CRCQ.

We say $\bar{x}$ satisfies the {\em MPODC-ERCPLD} if conditions (i)-(ii) are satisfied.  
\begin{itemize}
\item[{\rm (i)}]
$\{\nabla h_i(x)\}_{i=1}^m$ has constant rank,  for any $x\in \mathbb{B}(\bar x)$,  a neighborhood of $\bar x$. 
\item[{\rm (ii)}] Suppose $\mathcal{I}_1\subseteq \{1,\cdots,m\}$ is the index set that $\{\nabla h_i(\bar x)\}_{i\in \mathcal{I}_1}$  is a basis of $ \mbox{span} \{\nabla h_i(\bar x)\}_{i=1}^m.$ 
Suppose   there exist $\mathcal{I}_2\subseteq \bar I_g$, $L\subseteq \{1,\cdots,l\}$, 
 $v_i \in \partial g_i(\bar x), i\in \mathcal{I}_2$, $(\bar \eta_i^1,\bar \eta_i^2)\in \mathcal{N}_{{C_1\cup C_2}}(G_i(\bar x),H_i(\bar x))$, $i\in L$ such that 
\begin{eqnarray*}
({\{v_i, i\in \mathcal{I}_2\},\ \left\{\nabla h_i(\bar{x}), i\in \mathcal{I}_1;\ \nabla G_i(\bar x), i\in \mathcal{I}_3;\ \nabla H_i(\bar x), i\in \mathcal{I}_4\right\}})
\end{eqnarray*}  
is positive linearly dependent.
Then for sufficiently large $k$,
\begin{eqnarray*}
\{v_i^k\}_{i\in \mathcal{I}_2}\cup\{\nabla h_i(x^k)\}_{i\in \mathcal{I}_1}\cup  \{\nabla G_i(x^k)\}_{i\in \mathcal{I}_3}\cup\{\nabla H_i(x^k)\}_{i\in \mathcal{I}_4}
\end{eqnarray*} 
is linearly dependent, for any 
$\{x^k\}, \{v^k\}$ {and $\mathcal{I}_3, \mathcal{I}_4$ are defined as in MPODC-CRCQ.}
\end{itemize}

We say $\bar{x}$ satisfies the {\em MPODC-CPLD}, 
if all conditions in MPODC-RCPLD (ii) holds with $\mathcal{I}_1 \subseteq \{1,\cdots,m\}$ being an arbitrary set.
\end{defn}

The MPODC-quasi normality was defined in Benko et al. \cite[Corollary 5.1]{bch} and shown to be sufficient for error bounds \cite[Theorem 3.5]{bch}.
It is easy to prove that  MPODC-CPLD implies the MPODC-quasi normality similarly as the proof of Ye and Zhang \cite[Theorem 2.3]{yz14} and thus  MPODC-CPLD implies the error bound property.  However we do not know whether or not MPDC-CPLD  implies the quasi normality for MPDC. So we leave it as an open question.
It is easy to obtain the relationships among these constraint qualifications as in Figure 4.

\begin{figure}%
\centering
	\scriptsize
		 \tikzstyle{format}=[rectangle,draw,thin,fill=white]
		
		 \tikzstyle{test}=[diamond,aspect=2,draw,thin]
		
		\tikzstyle{point}=[coordinate,on grid,]
\begin{tikzpicture}
        \node[format](MPODC-PRCPLD){MPODC-PRCPLD};
        \node[format,right of=MPODC-PRCPLD,node distance=30mm](error bound){error bound};
         \node[format,right of=error bound,node distance=33mm](MPODC-quasi normality){MPODC-quasi normality}; 
        \node[format, below of=MPODC-PRCPLD,node distance=10mm](MPODC-NNAMCQ){MPODC-NNAMCQ};
        
        \node[format,left of=MPODC-NNAMCQ,node distance=30mm](MPODC-LICQ){MPODC-LICQ};
        \node[format,right of=MPODC-NNAMCQ,node distance=30mm](MPODC-CPLD){MPODC-CPLD};
	\node[format,right of=MPODC-CPLD,node distance=30mm](MPODC-RCPLD){MPODC-RCPLD};
       \node[format,right of=MPODC-RCPLD,node distance=27mm](M-stationarity){M-stationarity};
       \node[format, below of=MPODC-NNAMCQ,node distance=10mm](MPODC-CRCQ){MPODC-CRCQ};
       \node[format,right of=MPODC-CRCQ,node distance=30mm](MPODC-RCRCQ){MPODC-RCRCQ};
       \node[format,right of=MPODC-RCRCQ,node distance=30mm](MPODC-ERCPLD){MPODC-ERCPLD};
\draw[->](MPODC-LICQ)--(MPODC-NNAMCQ);
\draw[->](MPODC-ERCPLD)--(MPODC-RCPLD);
\draw[->](MPODC-PRCPLD)--(MPODC-RCPLD);
\draw[->](MPODC-PRCPLD)--(error bound);
\draw[->](MPODC-LICQ)--(MPODC-PRCPLD);
\draw[->](MPODC-NNAMCQ)--(MPODC-CPLD);
\draw[->](MPODC-CPLD)--(MPODC-RCPLD);
\draw[->](MPODC-RCPLD)--(M-stationarity);
\draw[->](MPODC-LICQ)--(MPODC-CRCQ);
\draw[->](MPODC-CRCQ)--(MPODC-RCRCQ);
\draw[->](MPODC-CRCQ)--(MPODC-CPLD);
\draw[->](MPODC-RCRCQ)--(MPODC-ERCPLD);
\draw[->](MPODC-CPLD)--(MPODC-quasi normality);
\draw[->](MPODC-quasi normality)--(error bound);
\end{tikzpicture}

\centering{Fig.4 Relation among constraint qualifications for MPODC}
\end{figure}

\subsection{RCPLD and piecewise RCPLD for MPECs}
In this subsection we apply our results for the ortho-disjunctive system to the equilibrium system:
\begin{eqnarray}
({\rm ES})~~~
&& g(x)\leq 0, h(x)=0, 
0\leq G(x)  \perp H(x) \geq 0 \label{feasibilitypo},
\end{eqnarray}
where the functions 
$ g:\mathbb{R}^d\to \mathbb{R}^n$ is locally Lipschitz continuous and $h:\mathbb{R}^d\to \mathbb{R}^m$, 
$G, H:\mathbb{R}^d\rightarrow \mathbb{R}^p$ are continuously differentiable.  For simplicity we have omitted the equality and inequality constraints. 

Because MFCQ fails at any feasible point of the MPEC if it is formulated as a nonlinear programming problem (see \cite[Proposition 1.1]{yzz}), various MPEC tailored  optimality conditions {were}
 investigated in the literature; see, e.g.,  \cite{y05}. 
Standard constraint qualifications,
such as LICQ, NNAMCQ, CRCQ, RCRCQ, CPLD and RCPLD,
 were extended to MPECs; see, e.g.,   \cite{gly1,hks13,y00,yz14}.
Moreover, the local error bound property was studied under the MPEC quasi-normality in \cite[Theorem 5.2]{gyz}, which is implied by the MPEC-CPLD.


Denote by ${\cal F}_{E}$ the  solution set of system  (\ref{feasibilitypo}). For  $\bar x\in {\cal F}_{E}$, let
\begin{eqnarray*}
&& {I}^{0+}:=\{i=1,\cdots,p: G_i(\bar x)=0, H_i(\bar x) >0\},\nonumber\\
&&{I}^{00}:=\{i=1,\cdots,p: G_i(\bar x)= H_i(\bar x)=0 \},\label{index}\\
&& {I}^{+0}:=\{i=1,\cdots,p: G_i(\bar x)>0, H_i(\bar x) =0 \}.\nonumber
\end{eqnarray*}
Around the point $\bar x$, the constraints $G_i(x)=0, i\in I^{0+}, H_i(x)=0 , i\in  I^{+0}$ can be treated as equality constraints and thus we can equivalently reformulate  the  system (\ref{feasibilitypo}) locally as 
\begin{eqnarray} \label{equv}
&&
g(x)\leq 0,\ h(x)=0,\nonumber\ 
G_i(x)=0, i\in I^{0+},\ H_i(x)=0 , i\in  I^{+0},\\
&& \Phi_i(x):=(G_i(x),H_i(x))\in\Omega_{E}, i\in I^{00}. 
\end{eqnarray}
Since $\Omega_{E}=(\mathbb{R}_+\times \{0\})\cup (\{0\}\times  {\mathbb{R}}_+)$, 
the equilibrium system  is an ortho-disjunctive system. 
From easily calculation, ${\cal N}_{\{0\}\times \mathbb{R}_{+}}(0,0)=\mathbb{R} \times \mathbb{R}_{-}$, ${\cal N}_{\mathbb{R}_{+}\times \{0\}}(0,0)=\mathbb{R}_{-}\times\mathbb{R}$.
Let $\mathcal{I}_3\subseteq \{i\in I^{00}:\eta_i^G\neq 0\}$ and $\mathcal{I}_4\subseteq \{i\in I^{00}:\eta_i^H\neq 0\}$.  By applying  MPODC-RCPLD and MPODC-PRCPLD to (\ref{equv}), 
we derived  the MPEC-RCPLD and MPEC-piecewise RCPLD as follows. Note that from Example \ref{exme}, for any $(\bar{\eta}_i^1,\bar{\eta}_i^2) \in \mathcal{N}_{\Omega_{E}}(0,0)$, $i\in I^{00}$, $\mathcal{I}_3$ contains the indexes $i\in I^{00}$ such that Cases (i) and (iii) happen, $\mathcal{I}_4$ contains the indexes $i\in I^{00}$ such that Cases (i) and (ii) happen.


\begin{defn}[MPEC-RCPLD and MPEC-PRCPLD] \label{rcpldmp}
We say the  {\em MPEC-RCPLD} holds at $\bar x\in {\cal F}_{E}$, if conditions (i)-(ii) are satisfied.
\begin{itemize}
\item[{\rm (i)}] The vectors $\{\nabla h_i(x)\}_{i=1}^m\cup \{\nabla G_i(x)\}_{i\in I^{0+}}\cup \{\nabla H_i(x)\}_{i\in I^{+0}}$ have the same rank for all $x$ in a small neighborhood of  $\bar x$.
\item[{\rm (ii)}] Let  $\mathcal{I}_1\subseteq \{1,\cdots,m\}$, $\mathcal{I}_2\subseteq I^{0+}$,  $\mathcal{I}_3\subseteq I^{+0}$ be such that  the set of vectors $\{\nabla h_i(\bar x)\}_{i\in \mathcal{I}_1}\cup \{\nabla G_i(\bar x)\}_{i\in \mathcal{I}_2}\cup \{\nabla H_i(\bar x)\}_{i\in\mathcal{I}_3}$ 
 is a basis for  
$$\mbox{ span } \{\{\nabla h_i(\bar x)\}_{i=1}^m\cup \{\nabla G_i(\bar x)\}_{i\in I^{0+}}\cup \{\nabla H_i(\bar x)\}_{i\in I^{+0}}\}.$$
Suppose there exist  index sets $\mathcal{I}_4\subseteq \bar{I}_g$, $\mathcal{I}_5,\mathcal{I}_6\subseteq I^{00}$, a nonzero vector $(\lambda^g,\lambda^h, \lambda^G,\lambda^H) \in \mathbb{R}^{n}\times\mathbb{R}^{m}\times \mathbb{R}^{p}\times\mathbb{R}^p$ satisfying $\lambda_i^g\geq 0$, 
$v_i \in \partial g_i(\bar x)$ for $i\in \mathcal{I}_4$ such that
\begin{eqnarray}
&&0 = \sum_{i\in \mathcal{I}_4} \lambda_i^g v_i +\sum_{i\in \mathcal{I}_1} \lambda_i^h \nabla h_i(\bar x) +\sum_{i\in\mathcal{I}_2\cup \mathcal{I}_5} \lambda_i^G \nabla G_i(\bar x) +\sum_{i\in\mathcal{I}_3\cup \mathcal{I}_6}\lambda_i^H \nabla H_i(\bar x),\label{mpec1}\\
&&\lambda_i^G< 0, \lambda_i^H< 0,\mbox{ or } \lambda_i^G \lambda_i^H=0,\ i\in I^{00}. \label{mpec2}
\end{eqnarray}
 \end{itemize}
Then the set of vectors 
\begin{eqnarray}
\{v_i^k\}_{i\in \mathcal{I}_4}\cup \{\nabla h_i(x^k)\}_{i\in \mathcal{I}_1}\cup \{\nabla G_i(x^k)\}_{i\in \mathcal{I}_2\cup\mathcal{I}_5}\cup\{\nabla H_i(x^k)\}_{i\in \mathcal{I}_3\cup\mathcal{I}_6}\label{ldvectors3}
\end{eqnarray}
 is linearly dependent for all sequences $\{x^k\},\{ v_i^k\}$ satisfying $x^k\rightarrow \bar x$, $x^k\neq \bar x$, $v_i^k \in \partial g_i(x^k)$, $v_i^k\rightarrow v_i$.

The  {\em MPEC-piecewise RCPLD (MPEC-PRCPLD)} holds at $\bar x\in {\cal F}_{E}$ if all conditions in MPEC-RCPLD hold with (\ref{mpec2}) replaced by
\begin{eqnarray}
&&\mbox{either } \lambda_i^G\leq 0\mbox{ or } \ \lambda_i^H\leq 0,\mbox{ for }\ i\in I^{00}. \label{mpec3}
\end{eqnarray}
\end{defn}

Note that when the nonsmooth inequality constraints are omitted,
MPEC-RCPLD  recovers  the one defined in Guo et. al \cite[Definition 3.4]{gly1}.
By applying our MPODC  constraint qualifications to MPECs, the MPEC-CRCQ (Hoheisel et al. \cite[Definition 2.4]{hks13}), MPEC-RCRCQ (Guo et al. \cite[Definition 3.4]{gly1}), MPEC-ERCPLD/MPEC-{$\tilde{r}$CPLD} (Guo et al. \cite[Definition 5.2]{gzl}, Chieu and Lee \cite[Definition 3.1]{cl14})
 and MPEC-CPLD can be derived.
 MPEC-CPLD we obtained is the same with the one defined in Ye and Zhang \cite[Definition 2.2(d)]{yz14}, which  is weaker than or equal to the original one given in Hoheisel et al. \cite[Definition 3.2]{hks12}.
Our MPEC-CPLD is neither implied or implies MPEC-{$\tilde{r}$CPLD}/MPEC-ERCPLD or MPEC-rRCPLD (Chieu and Lee \cite[Definition 3.1]{cl13}). 
But our MPEC-PRCPLD is weaker than or equal to MPEC-ERCPLD since no sign condition (\ref{mpec3}) is required for  {MPEC-ECPLD} to hold.
We summarize the relationship among various related constraint qualifications in Figure 4.

Guo et al. \cite[Theorem 5.1]{gzl} {showed} that  MPEC-RCPLD together with the strict complementarity condition implies the the error bound property. 
In this paper, by virtue of  Theorem \ref{xierr}, we have shown that  MPEC-PRCPLD is sufficient for the local error bound.  But MPEC-PRCPLD is stronger than MPEC-RCPLD. Hence  to our knowledge, 
the question whether or not MPEC-RCPLD (without any extra assumption such as the strict complementarity condition) would be sufficient for the error bound property is still open.

\begin{figure}%
\centering
\scriptsize
 \tikzstyle{format}=[rectangle,draw,thin,fill=white]
		
		 \tikzstyle{test}=[diamond,aspect=2,draw,thin]
		
		\tikzstyle{point}=[coordinate,on grid,]
\begin{tikzpicture}

        \node[format](MPEC-LICQ){MPEC-LICQ};
        \node[format,right of=MPEC-LICQ,node distance=28mm](MPEC-NNAMCQ){MPEC-NNAMCQ};
        \node[format,right of=MPEC-NNAMCQ,node distance=27mm](MPEC-CPLD){MPEC-CPLD};  
          \node[format,right of=MPEC-CPLD,node distance=26mm](MPEC-RCPLD){MPEC-RCPLD}; 
          \node[format,right of=MPEC-RCPLD,node distance=27mm](M-stationarity){M-stationarity};
          \node[format, below of=MPEC-LICQ,node distance=13mm](MPEC-CRCQ){MPEC-CRCQ};
       \node[format,right of=MPEC-CRCQ,node distance=26mm](MPEC-rCPLD){MPEC-rCPLD};       
        \node[format,right of=MPEC-rCPLD,node distance=27mm](MPEC-ERCPLD){MPEC-ERCPLD};
      \node[format,right of=MPEC-ERCPLD,node distance=28mm](MPEC-PRCPLD){MPEC-PRCPLD};
       \node[format,right of=MPEC-PRCPLD,node distance=27mm](error bound){error bound}; 
      
         \node[format, below of=MPEC-CRCQ,node distance=13mm](MPEC-RCRCQ){MPEC-RCRCQ};

\draw[->](MPEC-LICQ)--(MPEC-NNAMCQ);
\draw[->](MPEC-NNAMCQ)--(MPEC-CPLD);
\draw[->](MPEC-CPLD)--(MPEC-RCPLD);
\draw[->](MPEC-RCPLD)--(M-stationarity);
\draw[->](MPEC-LICQ)--(MPEC-CRCQ);
\draw[->](MPEC-CRCQ)--(MPEC-rCPLD);
\draw[->](MPEC-CRCQ)--(MPEC-CPLD);
\draw[->](MPEC-rCPLD)--(MPEC-ERCPLD);
\draw[->](MPEC-RCRCQ)--(MPEC-ERCPLD);
\draw[->](MPEC-ERCPLD)--(MPEC-PRCPLD);
\draw[->](MPEC-PRCPLD)--(MPEC-RCPLD);
\draw[->](MPEC-CRCQ)--(MPEC-RCRCQ);
\draw[->](MPEC-PRCPLD)--(error bound);
\draw[->](MPEC-CPLD)--(error bound);
\draw[->](error bound)--(M-stationarity);
\end{tikzpicture}

\centering{Fig.4 Relation among constraint qualifications for MPEC}
\end{figure}

\subsection{RCPLD and piecewise RCPLD for MPVCs}
In this subsection, we consider the vanishing {system}:
\begin{eqnarray*}
({\rm VS})~~~~~~~~ 
&& g(x)\leq 0, h(x)=0,\\
&& H_i(x)\geq 0, \quad G_i(x) H_i(x)\leq 0,\ i=1,\cdots,p,
\end{eqnarray*}
where 
$ g:\mathbb{R}^d\to \mathbb{R}^n,  i=1,\cdots,n$ are  assumed to be locally Lipschitz continuous and $h:\mathbb{R}^d\to \mathbb{R}^m, $
$G, H:\mathbb{R}^d\rightarrow \mathbb{R}^p$ are continuously differentiable.  For simplicity we have omitted the equality and inequality constraints. 

Denote by ${\cal F}_{V}$ the  solution set of system  (VS).
For $\bar x\in {\cal F}_{V}$, let
\begin{eqnarray*}
&& I^{+0}:=\{i: G_i(\bar x)>0, H_i(\bar x) =0\},\quad
 I^{00}:=\{i: G_i(\bar x)=0, H_i(\bar x) =0\},\\
&& I^{0+}:=\{i: G_i(\bar x)=0, H_i(\bar x) >0\},\quad
 I^{-0}:=\{i: G_i(\bar x)<0, H_i(\bar x) =0\},\\
&& I^{-+}:=\{i: G_i(\bar x)<0, H_i(\bar x) >0\}.
\end{eqnarray*}
 Around the point $\bar x$, we can equivalently formulate the constraint system of MPVC as 
\begin{eqnarray*}
&& g(x)\leq 0, h(x)=0, 
H_i(x)=0, i\in I^{+0}, G_i(x)\leq 0, i\in I^{0+}, H_i(x) \geq 0,  i\in I^{-0}, \nonumber\\
&& \Phi_i(x):=(G_i(x), H_i(x)) \in {\Omega}_{V}, i\in I^{00}. \label{VCsystem}
 \end{eqnarray*}

By applying the constraint qualifications for ortho-disjunctive system to the vanishing system, we derive the following conditions.

\begin{defn}\label{mvnrcpld}
We say $\bar x\in {\cal F}_{V}$ satisfies {\em MPVC-RCPLD}, if conditions (i)-(ii) are satisfied.
\begin{itemize}
\item[{\rm (i)}]
$\{\nabla h_i(x)\}_{i=1}^m\cup \{\nabla H_i(x)\}_{i\in I^{+0}}$ has constant rank, for any $x\in \mathbb{B}(\bar x)$,  a neighborhood of $\bar x$. 
\item[{\rm (ii)}]   Let  $\mathcal{I}_1\subseteq \{1,\cdots,m\}$, $\mathcal{I}_2\subseteq I^{+0}$ be such that the set of vectors $\{\nabla h_i(\bar x)\}_{i\in \mathcal{I}_1}\cup \{\nabla H_i(\bar x)\}_{i\in \mathcal{I}_2}$ 
 is a basis for  
$\mbox{ span}\{\{\nabla h_i(\bar x)\}_{i=1}^m\cup \{\nabla H_i(\bar x)\}_{i\in I^{+0}}\}.$
Suppose there exist index sets $\mathcal{I}_3\subseteq \bar{I}_g$,   $\mathcal{I}_4\in I^{0+}$, $\mathcal{I}_5\in I^{-0}$, $\mathcal{I}_6,\mathcal{I}_7\in I^{00}$, a nonzero vector $(\lambda^g,\lambda^h, \lambda^G,\lambda^H) \in \mathbb{R}^{n}\times\mathbb{R}^{m}\times \mathbb{R}^{p}\times\mathbb{R}^p$ and $v_i \in \partial g_i(\bar x)$,  for each $i\in \mathcal{I}_3$,
 such that 
 \begin{eqnarray}
&&0 = \sum_{i\in \mathcal{I}_3} \lambda_i^g v_i +\sum_{i\in \mathcal{I}_1} \lambda_i^h \nabla h_i(\bar x) +\sum_{i\in \mathcal{I}_4\cup \mathcal{I}_6}\lambda_i^G \nabla G_i(\bar x)+\sum_{i\in\mathcal{I}_2\cup \mathcal{I}_5\cup \mathcal{I}_7} \lambda_i^H \nabla H_i(\bar x),\nonumber \\
&&
\lambda_i^g\geq 0,\ i\in \mathcal{I}_3,\quad
\lambda_i^G\geq 0,\ \forall i \in I^{0+},\quad
\lambda_i^H\leq 0, 
\forall i \in I^{-0},\nonumber \\
&& \lambda_i^G \lambda_i^H=0,\ \lambda_i^G\geq 0,\ \forall i \in I^{00}.
\label{mpvc3}
\end{eqnarray}
\end{itemize}
Then the set of vectors 
$$\{v_i^k\}_{i\in \mathcal{I}_3}\cup \{\nabla h_i(x^k)\}_{i\in \mathcal{I}_1}\cup\{\nabla G_i(x^k)\}_{i\in \mathcal{I}_4\cup \mathcal{I}_6}\cup\{\nabla H_i(x^k)\}_{i\in \mathcal{I}_2\cup \mathcal{I}_5\cup \mathcal{I}_7},$$ 
where  $v_i^k \in \partial g_i(x^k)$, 
 is linearly dependent for all $\{x^k\}, \{v^k\}$ satisfying $x^k\rightarrow \bar x$, $x^k\neq \bar x$,  $v^k\rightarrow v$.
  
We say $\bar x\in {\cal F}_{V}$ satisfies  {\em MPVC-piecewise RCPLD (MPVC-PRCPLD)} if all conditions in MPVC-RCPLD hold with (\ref{mpvc3}) replaced by
\begin{eqnarray*}
{ \lambda_i^G=0 \mbox{ or }\ \lambda_i^G> 0,\ \lambda_i^H\leq 0, \mbox{ for }i\in I^{00}}.\label{mpvc4}
\end{eqnarray*}

\end{defn}
The MPVC-RCPLD is a constraint qualification for the {M-stationarity} from Theorem \ref{OC2}.
To our knowledge, MPVC-RCPLD has never been introduced before.

The error bound property for the MPVCs holds under the MPVC-piecewise RCPLD condition from Theorem \ref{xierr}.

\subsection{RCPLD for MPSCs}
In this subsection we apply our results to  the switching system:
\begin{eqnarray*}
({\rm SS})~~~~~~~~
&& g(x)\leq 0, h(x)=0, 
G_i(x) H_i(x)=0,\ i=1,\cdots,p,
\end{eqnarray*}
where 
$ g:\mathbb{R}^d\to \mathbb{R}^n,  i=1,\cdots,n$ are  assumed to be locally Lipschitz continuous and $h:\mathbb{R}^d\to \mathbb{R}^m$,  
$G, H:\mathbb{R}^d\rightarrow \mathbb{R}^p$ are continuously differentiable.
For simplicity we have omitted the equality and inequality constraints. 


Denote by ${\cal F}_{S}$ the  solution set of system  (SS).
For $\bar x\in {\cal F}_{S}$,  let
 \begin{eqnarray*}
 && I^{G}:=\{i: G_i(\bar x)=0, H_i(\bar x) \neq 0\},\ I^{H}:=\{i: G_i(\bar x)\neq 0, H_i(\bar x) =0\},\\
  && I^{GH}:=\{i: G_i(\bar x)=0, H_i(\bar x) =0\}.
\end{eqnarray*}
 Around $\bar x$, 
the switching system can be equivalently formulated as
\begin{equation}
g(x)\leq 0, h(x)=0,
G_i(x)=0, i\in I^{G}, H_i(x)=0 , i\in  I^{H}, (G_i(x),H_i(x))\in {\Omega}_{S}, i\in I^{GH}.\label{switchings}
\end{equation}

Let ${\Omega}_{S}=(\mathbb{R}\times \{0\})\cup (\{0\}\times\mathbb{R})$,
$
\displaystyle  \mathcal{N}_{\mathbb{R}\times \{0\}}(x_1,0)={\rm span}\{e_2\} 
$
 and 
$
\displaystyle\mathcal{N}_{\{0\}\times\mathbb{R}}(0,x_2)={\rm span}\{e_1\} 
$.
From (\ref{regularcone})-(\ref{limitingcone}),   the limiting normal cone 
{
\begin{eqnarray*}
\mathcal{N}_{\Omega_{S}}(0,0) 
&=& \bigcup_{x\in \mathbb{R}^2} \widehat{ \mathcal{N}}_{\Omega_{S}}(x)
 =\{ (\eta_1,\eta_2): \eta_1 \eta_2=0\}
 =\mathcal{N}_{C_1}(0,0)\cup \mathcal{N}_{C_2}(0,0).
\end{eqnarray*}
It is easy to versify that the MPSC-RCPLD coincides with the MPSC piecewise RCPLD.}

For any $( \eta_i^G,\eta_i^H) \in \mathcal{N}_{\Omega_{S}}(0,0)$, {$i\in I^{GH}$,} $ \eta_i^G,\eta_i^H$ can not be nonzero at the same time.
Let $\mathcal{I}_3\subseteq \{i\in I^{GH}:\eta_i^G\neq 0\}$ and $\mathcal{I}_4\subseteq \{i\in I^{GH}:\eta_i^H\neq 0\}$, then $\mathcal{I}_3\cap \mathcal{I}_4=\emptyset$.
Hence using these calculations and applying the constraint qualifications for ortho-disjunctive system to the system (\ref{switchings}), we obtain the following equivalent {definitions}. 

\begin{defn}\label{mpscnrcpld} 
We say $\bar x\in {\cal F}_{S}$ satisfies {\em MPSC-RCPLD}, if conditions (i)-(ii) are satisfied.
\begin{itemize}
\item[{\rm (i)}]
The vectors $\{\nabla h_i(x)\}_{i=1}^m\cup \{\nabla G_i(x)\}_{i\in I^{G}}\cup \{\nabla H_i(x)\}_{i\in I^{H}}$
has constant rank, for any $x\in \mathbb{B}(\bar x)$,  a neighborhood of $\bar x$. 
 \item[{\rm (ii)}]   Let  $\mathcal{I}_1\subseteq \{1,\cdots,m\}$, $\mathcal{I}_2\subseteq I^{G}$,  $\mathcal{I}_3\subseteq I^{H}$ be such that  the set of vectors $\{\nabla h_i(\bar x)\}_{i\in \mathcal{I}_1}\cup \{\nabla G_i(\bar x)\}_{i\in \mathcal{I}_2}\cup \{\nabla H_i(\bar x)\}_{i\in\mathcal{I}_3}$ 
 is a basis for  
$\mbox{ span } \{\{\nabla h_i(\bar x)\}_{i=1}^m\cup \{\nabla G_i(\bar x)\}_{i\in I^{G}}\cup \{\nabla H_i(\bar x)\}_{i\in I^{H}}\}.$
Suppose there exist index sets $\mathcal{I}_4\subseteq \bar{I}_g$,  $\mathcal{I}_5,\mathcal{I}_6\subseteq I^{GH}$ and a nonzero vector $(\lambda^g,\lambda^h, \lambda^G,\lambda^H) \in \mathbb{R}^{n}\times\mathbb{R}^{m}\times \mathbb{R}^{p}\times\mathbb{R}^p$ satisfying  $\lambda_i^g\geq 0,  v_i \in \partial g_i(\bar x),\ i\in \mathcal{I}_4$ and
\begin{eqnarray}
&&0 = \sum_{i\in \mathcal{I}_4} \lambda_i^g v_i +\sum_{i\in \mathcal{I}_1} \lambda_i^h \nabla h_i(\bar x) +\sum_{i\in\mathcal{I}_2\cup \mathcal{I}_5} \lambda_i^G \nabla G_i(\bar x) +\sum_{i\in\mathcal{I}_3\cup \mathcal{I}_6}\lambda_i^H \nabla H_i(\bar x),\nonumber \\
&& \lambda_i^G \lambda_i^H=0, \forall i \in I^{GH}.\label{mpsc2}
\end{eqnarray}
 \end{itemize}
Then the set of vectors 
$$\{v_i^k\}_{i\in \mathcal{I}_4}\cup \{\nabla h_i(x^k)\}_{i\in \mathcal{I}_1}\cup \{\nabla G_i(x^k)\}_{i\in \mathcal{I}_2\cup\mathcal{I}_5}\cup\{\nabla H_i(x^k)\}_{i\in \mathcal{I}_3\cup\mathcal{I}_6},$$ 
where  $v_i^k \in \partial g_i(x^k)$, 
 is linearly dependent for all $\{x^k\}, \{v^k\}$ satisfying $x^k\rightarrow \bar x$, $x^k\neq \bar x$,  $v^k\rightarrow v$.

\end{defn}

Since the MPSC-piecewise RCPLD coincides with the MPSC-RCPLD, from Theorem \ref{xierr}, the local error bound holds under MPSC-RCPLD.

\section{{Conclusions}}
In this paper, by using the normal cone structure of the disjunctive set, we have proposed RCPLD  for the disjunctive system which is 
{weaker than} the one proposed for the case where $\Gamma_i$ is an arbitrary closed set and is still a constraint qualification for the M-{stationarity} condition. Moreover we propose sufficient conditions for RCPLD  and study the relationships among them. 
Our new results  not only provide {a} uniform framework of the constraint qualifications for different disjunctive systems,
but {can} also derive weaker conditions for certain ortho-disjunctive systems than in the literature, such as MPVC-CRCQ, MPVC-CPLD and MPSC-CRCQ, MPSC-RCRCQ.
Moreover, we {can} also obtain new MPVC variants of RCRCQ, ERCPLD and RCPLD conditions that have never been defined before.
To obtain the local error bounds, we have proposed the piecewise RCPLD for the disjunctive system. 
{In particular, the  piecewise RCPLD derives the error bounds that does not need  the strict complementarity condition for MPECs.}





\end{document}